\documentclass{amsproc}

\usepackage{amssymb}
\usepackage{amscd}
\usepackage{amsthm}
\usepackage{graphicx}

\def\Zee{\mathbb{Z}}
\def\Q{\mathbb{Q}}
\def\Ar{\mathbb{R}}
\def\Cee{\mathbb{C}}

\def\Pee{\mathbb{P}}
\def\Rrr{\mathbb{R}}

\def\scrA{\mathcal{A}}
\def\scrB{\mathcal{B}}
\def\scrC{\mathcal{C}}
\def\scrE{\mathcal{E}}
\def\scrF{\mathcal{F}}
\def\scrG{\mathcal{G}}
\def\scrO{\mathcal{O}}
\def\scrS{\mathcal{S}}
\def\scrT{\mathcal{T}}
\def\scrZ{\mathcal{Z}}
\def\pa{\partial}
\def\tY{\widetilde{Y}}
\def\tX{\widetilde{X}}
\def\tW{\widetilde{W}}

\def\Ker{\operatorname{Ker}}

\def\Image{\operatorname{Image}}
\def\Tor{\operatorname{Tor}}

\def\Hom{\operatorname{Hom}}

\def\Id{\operatorname{Id}}

\def\int{\operatorname{int}}
\def\rk{\operatorname{rank}}
\def\codim{\operatorname{codim}}

\def\Spec{\operatorname{Spec}}
\def\Span{\operatorname{Span}}
\def\Rspan{\operatorname{{\mathbb{R}}-span}}
\def\Zspan{\operatorname{{\mathbb{Z}}-span}}

\newtheorem{theorem}{Theorem}[section]
\newtheorem{lemma}[theorem]{Lemma}
\newtheorem{corollary}[theorem]{Corollary}
\newtheorem{prop}[theorem]{Proposition}

\newtheorem{claim}[theorem]{Claim}
\theoremstyle{definition}
\newtheorem{definition}[theorem]{Definition}
\newtheorem{example}[theorem]{Example}

\theoremstyle{remark}
\newtheorem{remark}[theorem]{Remark}

\numberwithin{equation}{section}

\begin{document}

\title[Algebraic Topology of Calabi-Yau Threefolds in Toric
Varieties]{Algebraic Topology of Calabi-Yau Threefolds \\ in Toric
Varieties}

\author{Charles F. Doran}
\address{Department of Mathematics, University of Washington, Seattle,
Washington 98195}
\email{doran@math.washington.edu}

\author{John W. Morgan}
\address{Department of Mathematics, Columbia University, New York, New York
10027}
\email{jm@math.columbia.edu}

%    General info
%\subjclass{Primary 14D07, 14J32; Secondary 14M25, 19L64}
%\date{January 1, 1994 and, in revised form, June 22, 1994.}


%\keywords{Differential geometry, algebraic geometry}

\begin{abstract}
We compute the integral homology (including torsion), the topological K-theory, and the
Hodge structure on cohomology of Calabi-Yau threefold hypersurfaces and complete intersections
in Gorenstein toric Fano varieties.  The methods are purely topological.
\end{abstract}

\maketitle

\setcounter{tocdepth}{3}
\tableofcontents

\newpage

One of the most fruitful sources of Calabi-Yau threefolds is
hypersurfaces, or more generally complete intersections, in toric
varieties.  This is especially true since there is a proposal for
the mirror of any such Calabi-Yau three-fold. Usually the toric
varieties associated to convex lattice polytopes are singular,
causing the Calabi-Yau threefolds in them also to be singular, so
that to get smooth Calabi-Yau threefolds we must resolve the
ambient singularities and take the preimage in the resolution of
the singular Calabi-Yau threefold. This can be done torically by
the combinatorial device of taking a triangulation of the boundary
of the convex lattice polytope defining the toric variety where
the vertices of the triangulation are exactly the lattice points
contained in the boundary of the polytope.  In general, there will
be many such triangulations of a given lattice polytope, leading
to different ambient resolutions producing different families of
Calabi-Yau threefolds associated with the original toric variety.

In spite of the existence of many different such resolutions of a
given singular object, there are lattice-theoretic formulas for
the Hodge numbers of these resolutions expressed in terms of the
lattice polytope and its polar. Thus, the Hodge numbers of all the
different  resolutions coming from different triangulations are
the same. The proofs of these combinatorial formulas rely on the
Griffiths-Dwork method of computing Hodge numbers using residues
of meromorphic differentials on the complement of the Calabi-Yau
threefolds in the toric ambient space.

In this paper we study the resolution process from a more
topological point of view. In the topological study of these
objects, one treats both complete intersections and hypersurfaces
simultaneously.  This approach allows us to establish refined
versions of the Hodge number counts described above.  By directly
considering the topology and algebraic geometry of the resolutions
we compare the cohomology of the resolution with the cohomology of
the singular hypersurface or complete intersection. From this we
are able to see several things. The image of $H^3$ of the singular
object in $H^3$ of the resolution is identified with the weight 3
quotient of $H^3$ of the singular object and hence is the same for
all resolutions. The rest of $H^3$ of the resolution is of Hodge
types $(2,1)$ and $(1,2)$. It is described as a sum of Hodge
structures, summed over the edges of the polytope. The Hodge
structure associated to an edge is the tensor product of the Hodge
structures on $H^1$ of a smooth curve, contained in the closure of
the two-dimensional toric orbit given by that edge, with a Hodge
structure of type $(1,1)$ on the free abelian group with basis the
set of lattice points interior to the edge. It follows that the
rational Hodge structure on $H^3$ is independent of the choice of
resolution. Also, under the identification of $H^{2,1}$ with the
tangent space to the moduli space of complex structures, the image
of $H^{2,1}$ of the singular object is identified with the tangent
space of so called ``polynomial deformations'', i.e., the subspace
of deformations obtained by varying the hypersurface or complete
intersection in the toric variety.  The non-polynomial
deformations are then accounted for by the curves associated to
the edges of the polytope and the lattice points in the interiors
of these edges. In this way, we recover the so called ``correction
term'' describing the dimension of the non-polynomial deformations
directly from the resolution description.

It also follows from the resolution description that every class
in $H^2$ of the resolution is Poincar\'{e} dual to a divisor, a
toric divisor. As a consequence, the Hodge type of $H^2$ of the
resolution is $(1,1)$. We are able to recover the combinatorial
formula for the rank of this group from an understanding of the
resolutions of singularities.

Besides allowing us to enhance the combinatorial counts of Hodge
numbers to results about Hodge structures and allowing us to see
directly in terms of subvarieties and other topological objects
the sources of these homology groups, there are other advantages
to our topological approach.  It permits us to establish results
over the integers, and these of course lead to results in
topological K-theory. Since K-theory is the repository for the
most refined conjectures about mirror symmetry, understanding of
the integral homology and cohomology is crucial for more refined
tests of mirror symmetry.\footnote{In \cite{BK}, Batyrev and Kreuzer 
use a different method to study the integral cohomology of Calabi-Yau 
threefold hypersurfaces in toric varieties and obtain related
results, including some computations suggesting a role played by torsion 
in cohomology in mirror symmetry.}  While this paper establishes all the
necessary topology results over $\Zee$ in order to examine the
mirror symmetry proposals (see, for example, \cite{DM}), we do not
explore these issues here; they will be taken up in another paper.

\section{Preliminaries and statements of results}

Let $N$  be a lattice of dimension $n$, $\Delta \subset N_\Rrr = N
\otimes_\Zee \Rrr$
 be a reflexive polytope, and $\Pee_\Delta$
 the associated toric variety. Denote by $V=V(\Delta)$
 the set of
vertices of $\Delta$. If $n=4$, denote by $Y$ the vanishing locus
of a generic section of $\scrO(D_\infty)$ where
$$D_\infty=\sum_{v\in V}\overline{\scrO_v},$$
 and
$\scrO_v$  denotes the codimension-one torus orbit in
$\Pee_\Delta$ associated to the vertex $v$.

Here are the standard facts
 about toric varieties determined by reflexive polytopes
 \cite[\S3.5]{CK}
\begin{lemma}
Suppose that $\Delta\subset N_\Rrr$ is a reflexive polytope. Let
$D_\infty\subset \Pee_\Delta$ be the  Weil divisor
$$D_\infty = \sum_{v \in V(\Delta)} \overline{\scrO}_v.$$
The complete variety $\Pee_\Delta$ is Cohen-Macaulay. The divisor
$D_\infty$ is a Cartier divisor in $\Pee_\Delta$. The line bundle
$\scrO(D_\infty)$ is very ample, and anti-canonical.
 \end{lemma}

If $n>4$ suppose further that we have a NEF partition
$V=V_1\coprod V_2\coprod\cdots \coprod V_{n-3}$.  In this case,
suppose
$$D_i = \sum_{v \in V_i} \overline{\scrO_v}$$
 are ample divisors in $\Pee_\Delta$,  and let $Y\subset \Pee_\Delta$  be the complete intersection of generic sections of
$\scrO(D_i)$, $1 \leq i \leq n-3$.

For both the case when $Y$ is a hypersurface and a complete
intersection, we denote by $Y_i \subset Y$  the subvariety which
is the intersection of $Y$ with the union of the torus orbits of
codimension $\leq i$.

For $\scrT$  a maximal triangulation of $\pa\Delta$, i.e., one
whose set of vertices is $N \cap \pa\Delta$, let $\Pee_\scrT$
 be the toric variety
associated to the fan given by the cones over the simplices of
$\scrT$ and let
$$\rho \colon \Pee_\scrT \rightarrow \Pee_\Delta$$
be the morphism of toric varieties associated to the inclusion of
this fan into the fan determined by $\Delta$,  and define $\tY =
\rho^{-1}(Y)$.

The variety $\tY \subset \Pee_\scrT$ is a (smooth) Calabi-Yau
three-fold, and this construction of Calabi-Yau threefolds as
complete intersections in Gorenstein toric Fano varieties, due to
Batyrev and Borisov, is by far the richest known source of compact
Calabi-Yau threefolds.  In the hypersurface case alone, where the
requisite combinatorial data consists simply of a four dimensional
reflexive polytope $\Delta$, the classification of such polytopes,
due to Kreuzer and Skarke  \cite{KS}, results in more than 400
million example Calabi-Yau threefold hypersurface families.

Our goal in this paper is a complete understanding of the integral
(co)homology, Hodge structure, and topological K-theory of such
Calabi-Yau threefolds $\tY$.  The approach taken is thoroughly
algebro-topological by looking directly at this resolution.

The Hodge diamonds of the Calabi-Yau threefolds $\tY$ look like
this
$$\begin{array}{ccccccc}
  &   &         &    1    &         &   &   \\
  &   & 0       &         & 0       &   &   \\
  & 0 &         & h^{1,1} &         & 0 &   \\
1 &   & h^{2,1} &         & h^{2,1} &   & 1 \\
  & 0 &         & h^{1,1} &         & 0 &   \\
  &   & 0       &         & 0       &   &   \\
  &   &         &    1    &         &   &
  \end{array}$$
and there are combinatorial formulas expressing the Hodge numbers
$h^{1,1}$ and $h^{2,1}$. According to \cite{Bat,BB}, 
these formulas involve the distribution of
the points of $\Delta \cap N$ in the various faces of $\Delta$ and
its polar polytope $\Delta^\circ$  and also, in the case $n>4$,
the NEF partition/dual NEF partition. Let $E$, $F$, and $G$ denote
respectively the sets of edges, two-faces, and three-faces of
$\pa\Delta$, and $E^\circ$, $F^\circ$, and $G^\circ$ the same for
$\pa\Delta^\circ$. For simplicity, let us consider formulas for
  $h^{1,1}$ and  $h^{2,1}$
in the case of Calabi-Yau threefold hypersurfaces:
\begin{equation}
h^{1,1} = \ell(\Delta) - 5 - \sum_{g \in G} \ell^*(g) + \sum_{f
\in F} \ell^*(f) \ell^*(f^\vee) \label{h11hypintro}
\end{equation}
\begin{equation}
h^{2,1} = \ell(\Delta^\circ) - 5 - \sum_{v \in V} \ell^*(v^\vee) +
\sum_{e \in E} \ell^*(e) \ell^*(e^\vee). \label{h21hypintro}
\end{equation}
Here, $\ell(\Delta)$  is the cardinality of $N\cap \Delta$, and
$\ell(\Delta^\circ)$ is defined analogously with respect to the
dual lattice $\Hom(N,\Zee)$.  For a face $\alpha$ of $\pa \Delta$,
$\ell^*(\alpha)$  denotes the number of lattice points in the
relative interior of $\alpha$, and $\alpha^\vee$
 is the dual face to $\alpha$ in
$\pa\Delta^\circ$.

Although these formulas were used by Batyrev to check predictions from mirror
symmetry on the level of Hodge diamonds, there are two difficulties with this
approach which prevent the broader applicability of its methods.
First, for many applications to physical questions one would
really like to have a direct geometric construction of the
topological cycles generating the cohomology, especially over
$\Zee$.  One would like a complete description of the Hodge
structure on cohomology, rather than just the Hodge numbers, and
to have it expressed directly in terms of these geometric
representatives. Also, there is the issue of relating these Hodge
structures for $\tY$ to those of the original singular variety
$Y$, and seeing how, if at all, the answer depends on the choice
of maximal triangulation $\scrT$.

We begin by describing our results in the case of $H^2(\tY)$,
where the geometric representatives are more similar to those in
the literature, though the algebro-topological methods of proof
are quite different.

As is well known, all of $H^2(\tY)$ is spanned by algebraic cycles
and hence the Hodge decomposition of $H^2(\tY)$ is all of type
$(1,1)$. We establish in fact
\begin{theorem} \label{thmH2intro}
Let $\Delta$ be a reflexive polytope and $\scrT$ a maximal
triangulation of $\partial \Delta$. If ${\rm dim}(\Delta)=4$, let
$Y\subset \Pee_\Delta$ be a generic section of $\scrO(D_\infty)$.
If ${\rm dim}\Delta>4$ suppose that we have a NEF partition of
$V(\Delta)$ with associated divisors $D_i$, with
 the divisors $D_i \subset \Pee_\Delta$ being ample.
Let $\rho\colon \Pee_\scrT\to \Pee_\Delta$ be the natural resolution and let
$\tY=\rho^{-1}(Y)$. We have the following results.
\begin{enumerate}
\item In the hypersurface case,  $$\rk H^2(\tY) = \sum_{e \in E}
\ell^*(e) + \sum_{f \in F} \ell^*(f) \cdot (\ell^*(f^\vee) + 1) +
\# V - 4 .$$ \item The Hodge structure on $H^2(\tY)$ is of type
$(1,1)$. Every integral class in $H^2(\tY; \Zee)$ is an integral
linear combination of classes Poincar\'e dual to irreducible
components of toric divisors in $\tY$, i.e., divisors in $\tY$
that are irreducible components of the intersection of $\tY$ with
the closure of a codimension-one torus orbit in $\Pee_{\scrT}$.
 \item The mixed Hodge structure
on $H^2(Y)$ is a pure Hodge structure of weight 2 and Hodge type
$(1,1)$.
\end{enumerate}
\end{theorem}

\begin{remark}
Let us describe the irreducible components of the intersection of
$\tY$ with the closures of a codimension-one torus-orbit. These
orbits are indexed by the vertices of $\scrT$, i.e., $N\cap
\partial \Delta$. For each $\ell \in \pa \Delta \cap N$, there is
a divisor $R(\ell) \subset \tY$.  Let $\widetilde{\scrO}_\ell$ be
the associated torus orbit in $\Pee_{\scrT}$. Then $R(\ell)$ is
the closure of the intersection $\tY \cap \widetilde{\scrO}_\ell$.

When $\ell \in V$ is a vertex, $R(\ell)$ is an irreducible complex
surface projecting birationally onto the surface  $Y\cap
\overline{\scrO_\ell}$ in $Y$.

 When
$\ell \in \, \stackrel{\circ}{e}$ lies in the interior of an edge
$e \in E$, letting $\widehat{Z}(e)$ denote the smooth and complete
curve $Y \cap \overline{\scrO_e}$, $R(\ell)$ is a ruled surface
over the curve $\widehat{Z}(e)$ with generic fiber a smooth
rational curve $\Pee^1(\ell)$. In particular,
$H_3(R(\ell))=H_1(\widehat{Z}(e))\otimes H_2(\Pee^1(\ell))$. The
union over all $\ell\in \stackrel{\circ}{e}$ of the $R(\ell)$
produces a fibration where the generic fiber is
$A_{\ell^*(e)}$-configuration of rational curves.

Finally, when $\ell \in \stackrel{\circ}{f}$ for a 2-face $f \in F$, $R(\ell)$
is a disjoint union of $\# (Y \cap \scrO_f)$ copies of a surface and maps to $Y
\cap \scrO_f$ in $Y$.
\end{remark}

The results in the hypersurface case described in part (1) of the
theorem imply an expression for $h^{1,1}(\tY)$.  This agrees with
that from formula (\ref{h11hypintro}), since
$$\ell(\Delta) = 1 + \# V + \sum_{e \in E} \ell^*(e) + \sum_{f \in F}
\ell^*(f) + \sum_{g \in G} \ell^*(g) .$$ There is an analogue of
part (1) in the complete intersection case, and corresponding
expressions for $h^{1,1}(\tY)$ in the complete intersection case
(see Section \ref{subsectfirst} below and \cite[\S 8]{BB}.

We now describe our corresponding results for $H^3(\tY)$:

\begin{theorem} \label{thmH3intro}
In the hypersurface case or in the complete intersection case
when all the divisors $D_a\subset \Pee_\Delta$ are ample we have the following
results.
\begin{enumerate}
\item The image of $$\rho_* \colon H_3(\tY) \rightarrow H_3(Y)$$
is equal to the image of
$$H_3(Y_2) \rightarrow H_3(Y) \ .$$
\item    Let $\scrA_{\ell^*(e)}$ be the lattice (with symmetric
bilinear pairing) associated to the root system $A_{\ell^*(e)}$.
Then
$$\Ker \rho_* \colon H_3(\tY) \to H_3(Y)$$
is identified with $$\bigoplus_{e \in E} H_1(\widehat{Z}(e)) \otimes
\scrA_{\ell^*(e)} .$$  Under this isomorphism, the restriction to $\Ker \rho_*$
of the usual homological intersection pairing is identified with the direct sum
of the tensor products of the intersection pairings on $H_1(\widehat{Z}(e))$
with the natural pairing on $\scrA_{\ell^*(e)}$.
\end{enumerate}
\end{theorem}

It follows that
$$\Image \left( H_3(\tY; \Zee) \rightarrow H_3(Y; \Zee) \right)$$
is independent of the choice of maximal triangulation $\scrT$.  By
general results in mixed Hodge theory \cite{Del} in fact
$$\Image \left(H_3(\tY; \Q) \rightarrow H_3(Y; \Q) \right)$$
is $W_{-3} \subset H_3(Y; \Q)$ under the mixed Hodge structure on $H_3(Y;\Q)$.

\begin{corollary} \label{corH3introseq}
We have the following exact sequence of $\Q$-Hodge structures of
weight three:
$$0 \to \frac{H^3(Y; \Q)}{W_2H^3(Y;\Q)} \to H^3(\tY;\Q) \to
\bigoplus_{e \in E} H^1(\widehat{Z}(e)) \otimes \scrA^*_{\ell^*(e)} \to 0 \ ,$$
where the first term is a pure weight three quotient of $H^3(Y;\Q)$.  The
sequence splits so that the Hodge structure on $H^3(\tY;\Q)$ is isomorphic to
the direct sum of this pure weight three Hodge structure and the direct sum of
tensor products of Hodge structures from the third term. The factor
$\scrA^*_{\ell^*(e)}$ is defined to be of type $(1,1)$, so that the third term
has Hodge types $(2,1)$ and $(1,2)$.
\end{corollary}
\begin{corollary} \label{corH3introres}
The isomorphism type of the $\Q$-Hodge structure on $H^3(\tY)$ is
independent of the choice of resolution $\tY$ of $Y$
(corresponding to the choice of maximal triangulation $\scrT$ of
$\pa\Delta$).
\end{corollary}

The last terms in formulas (\ref{h11hypintro}) and
(\ref{h21hypintro}) are typically called the ``correction terms'',
and without these terms the formulas describe the so-called {\em
toric part} $h^{1,1}_{toric}$  of $h^{1,1}(\tY)$ and the {\em
polynomial part} $h^{2,1}_{poly}$ of $h^{2,1}(\tY)$. The number
$h^{1,1}_{toric}$ is the dimension of the subspace of
$H^{1,1}(\tY; \Cee)$ generated by the restriction to $\tY$ of the
$T_N$-invariant divisors on $\Pee_\scrT$.
%By part (3) of
%Theorem~\ref{thmH2intro}, $\tY \cap \widetilde{\scrO}_\ell$
%consists of a single irreducible component when $\ell$ is a vertex
%or $\ell \in \stackrel{\circ}{e}$ for $e \in E$, and when $\ell
%\in \stackrel{\circ}{f}$, we get a number of components equal to
%$\# \tY \cap \scrO_f$.
Thus, by parts (2) and (3) of Theorem~\ref{thmH2intro}, the
difference between $h^{1,1}_{toric}$ and $h^{1,1}$ is given by
$\ell^*(f) \cdot \ell^*(f^\vee)$ and reflects the fact that, for
$\ell \in \, \stackrel{\circ}{f}$, $R(\ell)$ has $\ell^*(f^\vee) +
1$ components. The number $h^{2,1}_{poly}$ counts the dimension of
the space of {\em polynomial deformations} of $\tY \subset
\Pee_\scrT$, i.e., the space of deformations of complex structure
of $\tY$ determined by the hypersurfaces in the anti-canonical
linear system $|-K_{\Pee_\scrT}|$, or equivalently deformations of
$\tY$ induced by taking the preimages under $\rho$ of deformations
of $Y\subset \Pee_\Delta$.

\begin{corollary}
 The sequence in
Corollary~\ref{corH3introseq} is compatible with the subspace of
polynomial deformations of $\tY$ in the sense that the subspace
$\bigl(H^3(Y)/W_2H^3(Y)\bigr)^{2,1}$  of $H^{2,1}(\tY)$ is the
tangent space to the space of polynomial deformations of $\tY$.
\end{corollary}

The methods used in establishing the above results on (co)homology
of the Calabi-Yau threefolds $\tY$ involve both ``standard'' local
decomposition techniques, guided by the toric structure on
$\Pee_\scrT$, and some combinatorial topology about the two-faces
of the reflexive polytope $\Delta$ itself.

The methods apply not only to integral cohomology but also to
$K$-theory. To see this we  use the direct link between
topological K-theory and integral homology (including torsion). A
careful analysis of the 7-dimensional stage $BSU^{(7)}$ of the
Postnikov tower for $BSU$ leads us to the following technical
result.
\begin{theorem}
Let $M$ be a closed, oriented 6-manifold.  The reduced even
K-group $\widetilde {K^0}(M)$ is isomorphic to
$$\{(c_1,c_2,c_3) \in H^2(M;\Zee) \oplus H^4(M;\Zee) \oplus
H^6(M;\Zee) \vert Sq^2 c_2 = [c_3]_2 +c_1c_2+c_1^3\}$$ where the
isomorphism is given by taking the 1st, 2nd, and 3rd Chern class.
The odd K-group $K^1(M) = \widetilde{K^0}(\Sigma M)$ is
$$H^1(M;\Zee) \times [\Sigma M, BSU^{(7)}]$$ and we have an exact
sequence
$$0 \to H^5(M;\Zee) \to [\Sigma M, BSU^{(7)}] \to H^3(M;\Zee) \to 0 \ ,$$
where the extension is the pullback of the extension
$$0 \to H^5(M;\Zee) \to A \to H^3(M;\Zee) \to 0$$
with extension class given by the $Sq^2$ map
$$H^3(M;\Zee/2\Zee) \to \frac{H^5(M;\Zee)}{2 H^5(M;\Zee)} =
H^5(M;\Zee/2\Zee) \ .$$
%coming from the fibration $$K(\Zee,6) \to
%BSU^{(7)} \to K(\Zee,4)$$ with $k$-invariant $\delta Sq^2 \iota_4
%\in H^7(K(\Zee,4);\Zee)$.
\end{theorem}

For Calabi-Yau threefolds, independent of whether or not the manifolds are
constructed as hypersurfaces or complete intersections in a toric variety, we
establish the following.

\begin{corollary}
Let $M$ be a Calabi-Yau threefold.  The topological K-groups of
$M$ are expressed in terms of the integral homology groups of $M$:
$$K^0(M) \cong \Zee \oplus H^2(M;\Zee) \oplus H^4(M;\Zee) \oplus
2 \cdot H^6(M;\Zee) \ ,$$ and
$$K^1(M) \cong H^1(M;\Zee) \oplus H^3(M;\Zee) \oplus H^5(M;\Zee) \ .$$
\end{corollary}

\vspace{.1in}

\section{Calabi-Yau manifolds in toric varieties}
\label{toricstuff}

This section lays out notations and basic results from toric
geometry which we will use in subsequent sections. Good references
for many of these are the books \cite{CK,Ful}.

Let $N \simeq {\mathbb{Z}}^n$ be a lattice, $M =
\Hom(N,{\mathbb{Z}})$ the dual lattice, $N_{\mathbb{R}} = N
\otimes_{\mathbb{Z}} {\mathbb{R}}$ the associated real vector
space, and $T_N = N \otimes_{\mathbb{Z}} {\mathbb{C}}^*$ the
associated complex torus.

\subsection{Toric varieties: The affine case}

Let $c$ be a rational polyhedral cone in $N$. We say that $c$ is
{\em strongly convex} if $c \cap (-c) = \{0\}$. The {\em dual
cone} $\check{c}$ of a strongly convex rational polyhedral cone
$c$ is
$$\check{c} = \{x \in M_{\mathbb{R}} \, | \, \langle x,c \rangle
\geq 0 \} \ .$$ Define the associated open affine toric variety by
$$U(c,N) = \Spec{\Cee}[M \cap \check{c}] \ .$$
The action $T_N \times U(c,N) \to U(c,N)$ is dual to the morphism
\begin{eqnarray*} \Cee[M \cap \check{c}] & \to &
\Cee[M \cap \check{c}]
\otimes_\Cee \Cee[M] \\
\chi^m & \mapsto & (\chi^m \otimes \chi^m)
\end{eqnarray*}
where the characters $\chi^m \colon T_N \to \Cee^*$ for $m \in M
\cap \check{c}$ generate the $\Cee$-algebra $\Cee[M \cap
\check{c}]$. Since $c$ is strongly convex, $M \cap \check{c}$
spans $M_\Rrr$ as a real vector space and hence this action has a
free $T_N$-orbit that is open and dense.

 The faces of $c$ are themselves
polyhedral cones and the set of faces is a poset under inclusion,
which we denote by $f \prec f'$. Consider the $1$-dimensional
faces $r \prec c$ of
 the cone $c$.
The {\em generators} of $c$, denoted as a set by $\scrG(c)$, are
the generators of the semigroups $r \cap N$ as $r$ runs through
all the $1$-dimensional faces of $c$. We denote by
$\stackrel{\circ}{c}$ the relative interior
$$\stackrel{\circ}{c} \ = \, c \setminus
\bigcup_{ \{f \prec \, c \, | \, f \neq c\} } f \ .$$ Let $f$ be a
face of $c$. We define $K_f\subset M\cap \check{c}$ as follows:
$$K_f=\{m\in M\cap \check{c} \, \bigl|\bigr. \, \langle m,f\rangle=0\}.$$
Corresponding to the face $f$ of $c$ there is a closed
$T_N$-invariant subset in $U(c,N)$ defined by
$$\scrZ_f=\bigcap_{m\in (M\cap \check{c}) \setminus K_f}\{\chi^m=0\}.$$
This subset is the closure in $U(c,N)$ of the $T_N$-orbit
$\scrO_f$ (sometimes denoted $\scrO_{f,N}$ if we need to specify
the lattice $N$ to avoid confusion). This orbit is the open subset
of $\scrZ_f$ consisting of all points $z\in \scrZ_f$ with the
property  that for every $m\in K_f$ $\chi^m(z)\not= 0$. The union
of the orbits $\scrO_f$, as $f$ ranges over all the faces of $c$,
is $U(c,N)$. The closure $\overline{\scrO}_f$ is the union of
$\scrO_{f'}$ as $f'$ ranges over all faces  of $c$ satisfying
$f\prec f'$. The open dense orbit corresponds to the face of $c$
given by the point $\{0\}$.  At the other extreme, there is
exactly one closed orbit; it is $\scrO_c$.  We denote by
$U^*(c,N)$ the complement of this orbit. For maximal dimensional
cones we have the following result about the closed orbit and also
about the structure of the toric variety $U(c,N)$.
\begin{lemma}\label{fixedpoint}
Suppose that $c$ spans $N_{\mathbb{R}}$, i.e. the dimension of $c$
is $n$. The orbit $\scrO_c$ associated to the cone $c$ itself is a
fixed point of the $T_N$-action.  If $c$ is simplicial, then
$U(c,N) \simeq {\mathbb{C}}^n / \Gamma$, where
$$\Gamma = N / \Zspan(\scrG(c)) \ .$$
\end{lemma}

Now let us consider the case when $c$ is not of maximal dimension.
Let $\Rspan(c) \subset N_{\mathbb{R}}$ be the linear span of the
cone $c$ as a vector subspace. The {\em dimension} of $c$ is
defined to be the dimension of $\Rspan(c)$. We define
 $N_c = N \cap \Rspan(c)$. Of course, $c$ is a polyhedral cone in
 $(N_c)_{\mathbb R}$. We denote by $V(c)$ the toric variety $U(c,N_c)$.
 According to Lemma~\ref{fixedpoint} there is a unique fixed point,
 denoted $0_c$ for the $T_{N_c}$-action on $V(c)$.

Dual to the inclusion $N_c \hookrightarrow N$ is the surjection $M
\to M_c$, where $M_c = \Hom(N_c, \Zee)$. Denote by $K_c$ the
kernel of $M \to M_c$.  We have in fact the (non-canonically
split) exact sequence of abelian semigroups
$$0 \to K_c \to M \cap \check{c} \to M_c \cap \check{c} \to 0 .$$
This gives rise to ring homomorphisms $$\Cee[K_c] \to \Cee[M\cap
\check{c}] \to \Cee[M_c \cap \check{c}] \ .$$ Since $K_c$ is
identified with $\Hom(N/N_c, \Zee)$, and since $N_c \subset N$ is
a direct summand, we have $T_{N/N_c} = T_N/T_{N_c}$. Thus, we  the
dual morphisms are
\begin{equation} \label{UcNcfibration}
V(c) \to U(c,N) \to T_N/T_{N_c} \ .
\end{equation}
The  map $U(c,N) \to T_N/T_{N_c}$ is a fibration and $V(c)$ is the
fiber over the identity element. Since $N_c$ is a direct summand
of $N$, there is a (non-canonical) splitting $T_N/T_{N_c}
\hookrightarrow T_N$. Restricting the action $T_N \times U(c,N)
\rightarrow U(c,N)$ to the image of the splitting determines an
isomorphism $$V(c) \times T_N/T_{N_c} \stackrel{\psi}{\to} U(c,N)
\ .$$ We have the product action of $T_{N_c} \times (T_N/T_{N_c})$
on $V(c) \times (T_N/T_{N_c})$. The splitting identifies this
product torus with $T_N$ and with this identification the
isomorphism $\psi$ is equivariant.  The closed orbit $\scrO_c$ in
$U(c,N)$ is the product of $0_c\in V(c)$ with $T_N/T_{N_c}$.
\begin{lemma} \label{lemmaUcNccone}
There is a $\Cee^*$-action on $V(c)$ whose the only fixed point is
$0_c$ and such that for any $x \in V(c)$, $\lim_{\lambda \to 0}
\lambda \cdot x = 0_c$.  Restricting to $\Rrr^+ \subset \Cee^*$
produces the natural contraction action of a real cone-structure
on $V(c)$ with cone point $0_c$.
\end{lemma}
\begin{proof}
Let $n \in N_c$ be a primitive vector contained in the interior of
the cone $c$.  Associated to this vector is an embedding
$\imath_n\colon \Cee^* \hookrightarrow N_c \otimes_\Zee \Cee^*$.
We claim that for any $x \in V(c)$,
$$\lim_{\lambda \to 0} \imath_n(\lambda) \cdot x$$
is the fixed point $0_c$. The reason is that if $u \in M_c
\setminus \{0\}$ is $\geq 0$ on the cone $c$, then $\langle u, n
\rangle > 0$ and hence
$$\lim_{\lambda \to 0} \chi^u(\imath_n(\lambda) \cdot x) = 0 .$$
This means that the orbits of the $\Rrr^+ \subset \Cee^*$ action
give a flow on $V(c)$ contracting it to the cone point.
\end{proof}

\subsection{Toric varieties: The general case}

A {\em fan} $\scrF$ in $N_\Rrr$ is a finite collection of strongly
convex rational polyhedral cones in $N_\Rrr$ such that
\begin{enumerate} \item If $c \in \scrF$, then every face of
$c$ is also in $\scrF$. \item If $c_1, c_2 \in \scrF$, then $c_1
\cap c_2$ is a face of each of $c_1$ and $c_2$.
\end{enumerate}
The {\em support} of a fan $\scrF$ is the set
$$|\scrF| = \bigcup_{c \in \scrF} c \subset N_\Rrr \ .$$

We obtain the {\em toric variety} $\Pee_\scrF$ from a fan $\scrF
\subset N_\Rrr$ by gluing together the affine $T_N$-spaces
$U(c,N)$ for all $c \in \scrF$, where $U(c_1,N)$ and $U(c_2,N)$
are glued together along $U(c_1 \cap c_2, N)$. (To avoid ambiguity
arising from choice of lattice, we sometimes denote the toric
variety associated with a fan $\scrF \subset N_\Rrr$ by
$\Pee_{\scrF,N}$.) The variety $\Pee_\scrF$ is complete if and
only if the fan $\scrF$ is {\em complete}, i.e., $|\scrF| =
N_\Rrr$. By the proposition on integral closure in
\cite[p.~29]{Ful}, irreducible toric varieties are normal, and
hence nonsingular through codimension two.  It follows from the
affine case that there is a one-to-one correspondence between the
$T_N$-orbits and cones of $\scrF$, with the adherence relation
being the opposite of the face relation; i.e., for cones $c_1$ and
$c_2$,
\begin{equation}\label{face}
\scrO_{c_1} \subset \overline{\scrO}_{c_2} \Leftrightarrow  c_2
\prec c_1 \ .\end{equation}

Given a strongly convex rational polyhedral cone $c$, we define
$\scrF_c$ to be the fan consisting of $c$ and all its faces. In
this case we have $\Pee_{\scrF_c} = U(c,N)$. Also considering $c$
as a polyhedral cone in $(N_c)_\Rrr$, let $\scrF_c(N_c)$ denote
the fan in $(N_c)_\Rrr$ consisting of $c$ and all its faces. We
have $\Pee_{\scrF_c(N_c)} = U(c,N_c)$. For a general fan $\scrF
\subset N_\Rrr$ we have $\Pee_\scrF = \bigcup_{c\in \scrF}
U(c,N)$.

One natural way that complete fans arise is from convex lattice
polytopes $\Delta\subset N_\Ar$.  For us a convex lattice polytope
contains the origin in its interior and has points of $N$ as its
vertices. We denote by $V(\Delta)$ its set of vertices, by
$E=E(\Delta)$ its set of edges and by $F=F(\Delta)$ its set of
two-faces. We define the fan $\scrF(\Delta)$
associated\footnote{Our associated fan is {\em not} what is
referred to in \cite{CK} as the normal fan of $\Delta$; in our
notation, the normal fan of $\Delta \subset N_\Rrr$ is a fan in
$M_\Rrr$.} with $\Delta$ consisting of cones over the faces of
$\partial \Delta$ (by convention we include as a face the empty
face whose cone is the origin), i.e., over a face $\gamma$
consider the cone
$$c(\gamma) = \{ \lambda (x - x') \, | \, x \in \Delta, x' \in \gamma, \lambda \geq
0 \} \subset N_\Rrr \ .$$  We define the associated toric variety
$$\Pee_\Delta = \Pee_{\scrF(\Delta)} \ .$$
For any face $\gamma$ of $\partial\Delta$ we denote by $\scrO_\gamma$ the
$T_N$-orbit $\scrO_{c(\gamma)}$, by $U(\gamma,N)$ the affine toric variety
$U(c(\gamma),N)$, by $N_\gamma$ the sublattice $N_{c(\gamma)}$, and by
$V(\gamma)$ the affine toric variety $V(c(\gamma))$.

A notation we shall use throughout is the following: If $\gamma$ is a face of a
lattice polytope then $\ell^*(\gamma)$ denotes the number of lattice points
contained in the relative interior of the face; i.e,
$$\ell^*(\gamma)=\#\left\{N \, \cap \stackrel{\circ}{\gamma}\right\}.$$

\subsubsection{Refinements and induced maps: The affine case}

We say that a fan $\scrF'$ {\em refines} the fan $\scrF$ if each
cone of $\scrF'$ is contained in a cone of $\scrF$ and that the
supports of $\scrF'$ and $\scrF$ are the same.

Let's begin with the case of refining a single cone. Fix a
strongly convex rational polyhedral cone $c$ and suppose that
$\scrF'_c$ is a fan refining $\scrF_c$.

\begin{lemma}\label{LemmaRho}
There is a $T_N$-invariant morphism $\rho \colon
\Pee_{\scrF'_c}\to U(c,N)$ which is a birational isomorphism.
\end{lemma}
\begin{proof}
Let $c'$ be a cone of $\scrF'_c$.  Since $c'\subset c$, it follows
that $\check{c} \subset \check{c}'$. The map $U(c',N)\to U(c,N)$
is dual to the inclusion
$$\Cee[M\cap \check{c}]\to \Cee[M\cap \check{c}'].$$
It is easy to see that these affine  maps are $T_N$-equivariant and
are fit together to define the map $\rho$ as given in the statement.
\end{proof}

We wish to describe the pre-image of the orbit $\scrO_{c}$ under
the map $\rho \colon \Pee_{\scrF'_c}\to U(c,N)$. Because of the
$T_N$-equivariance, this pre-image is a union of $T_N$-orbits. In
order to describe this map, let us introduce the {\em dual cell
decomposition} $\Sigma(\scrF'_c)$ of the refinement. There is one
cell $\sigma'$ dual to each cone $c'$ of $\Sigma(\scrF'_c)$
meeting $\stackrel{\circ}{c}$. The dimension of $\sigma'$ is
$n-{\rm dim}(c')$. The cell $\sigma'$ dual to $c'$ is a face of
the cell $\sigma''$ dual to $c''$ if and only if $c''\prec c'$.
Notice that this cell complex is of dimension at most $n-1$. It
need not be of this dimension, nor even be of homogeneous
dimension. We denote by $c_\sigma$ the cone in $\scrF'_c$ dual to
the cell $\sigma\in \Sigma(\scrF'_c)$.

\begin{lemma}\label{preimage}
The pre-image $\rho^{-1}(\scrO_{c})$ is given as
$$\rho^{-1}(\scrO_{c}) = \bigcup_{\sigma \in \Sigma(\scrF'_c)}
\scrO_{c_\sigma}.$$ This bijective correspondence between the
$T_N$-orbits in $\rho^{-1}(\scrO_{c})$ and the cells of
$\Sigma(\scrF'_c)$ identifies the face relationship in
$\Sigma(\scrF'_c)$ with adherence relationship of orbits.
\end{lemma}
\begin{proof}
Let $c'$ be a cone in $\scrF'_c$. Then
$\rho(\scrO_{c'})=\scrO_{c}$ if and only if every $m\in
M\cap\check{c}$ that is not identically zero on $c$ is not
identically zero on $c'$. Since $c'$ contains a relative interior
point of $c$ and since $m\ge 0$ on $c$, this is clear. The fact
that $\scrO_{c_1}$ is contained in the closure of $\scrO_{c_2}$ if
and only if $(c_1)_\sigma$ is a face of $(c_2)_\sigma$ is clear
from Equation~(\ref{face}).
\end{proof}

%We also denote by $\rho$ the natural morphism from the toric
%variety $\Pee_{\scrF'_c(N_c)}$ (associated with the refinement of
%$\scrF_c(N_c)$) to $U(c,N_c)$.

\begin{figure}[ht]
\centering
\includegraphics[width=3in]{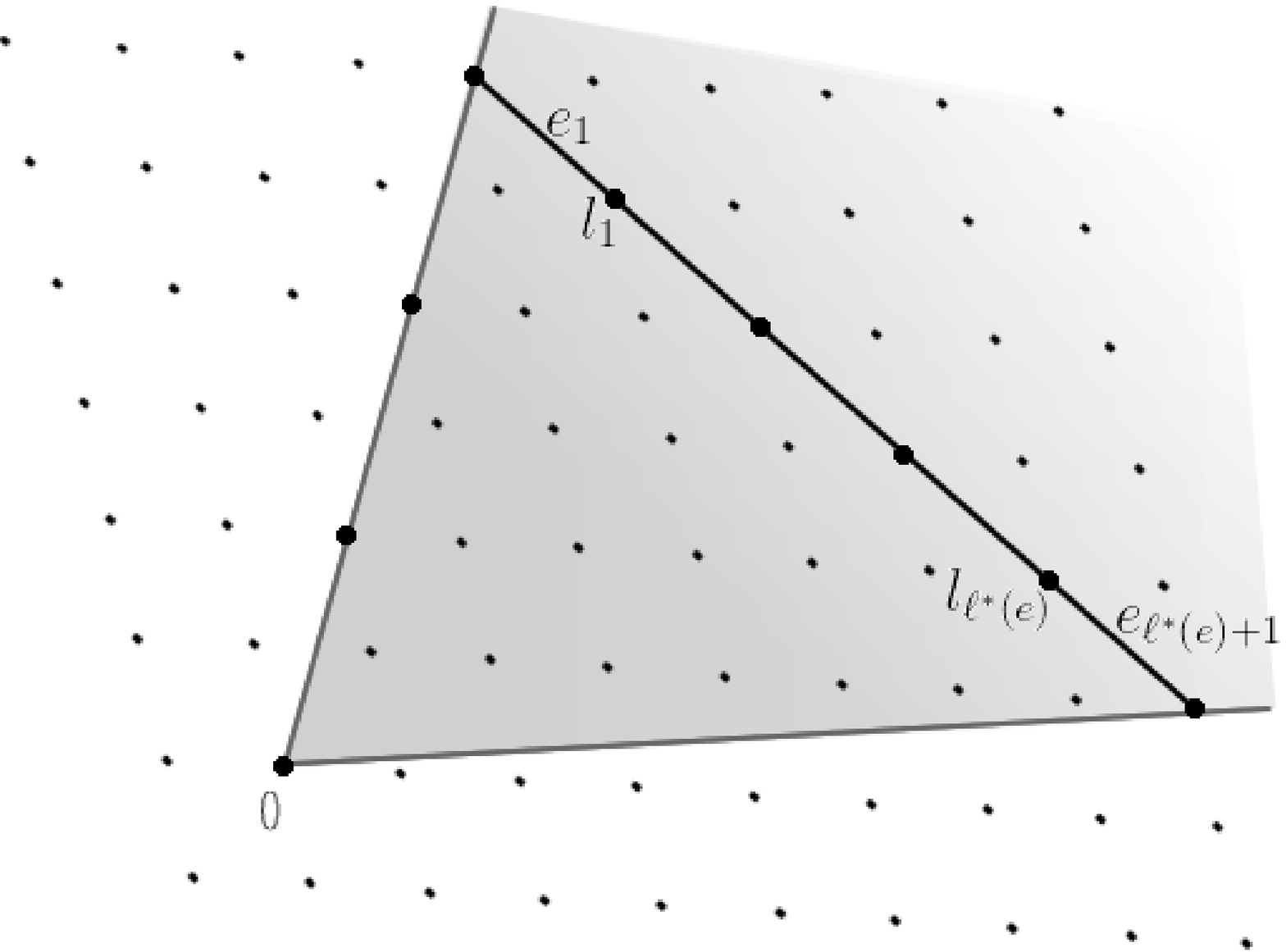}
\caption{Example~\ref{eg:edge}}
\end{figure}

\begin{example} \label{eg:edge}
%Let $\Delta$ be a reflexive polytope in $N_\Rrr$ and $e$ be an
%edge of $\Delta$.
Let $c$ be the cone on an edge $e$ with vertices indecomposable elements of
$N$. Set $k(e)=\ell^*(e)+1$. Then $V(e) \simeq {\Cee}^2 / (\Zee / k(e) \Zee)$.
The action is given by the primitive $k(e)^{th}$ root of unity $\zeta$ acting
by $\zeta(z_1,z_2) = (\zeta z_1, \zeta^{-1} z_2)$. In particular, the space
$V(e)$ is homeomorphic to $c(L(k(e), -1))$, the cone over the Lens space
$L(k(e), -1)$. Consider the subdivision $\scrF_e$ of $e$ given by the
decomposition
$$e = e_1 \cup \ldots \cup e_{\ell^*(e)+1}$$ with $l_i = e_i \cap e_{i+1}$ for
$i = 1, \ldots, \ell^*(e)$ the points $N \cap
\stackrel{\circ}{e}$. Then $\Pee_{\scrF'_e}$ is smooth and the
morphism $\rho\colon\Pee_{\scrF'_e} \to V(e)$ resolves the
singularity at $0_{e}$ into an $A_{\ell^*(e)}$ configuration of
rational curves --- a chain of $\ell^*(e)$ rational curves.
\end{example}

\begin{lemma}
Suppose that the $\Rrr$-span of $c$ is $N$, so that we have the
unique fixed point $0_c$ of the $T_N$-action on $U(c,N)$.  Then
the fiber $\rho^{-1}(0_c)$ is a compact variety, and
$\Pee_{\scrF'_c}$ deformation retracts to $\rho^{-1}(0_c)$.
%$\Pee_{\scrF'_c}$ deformation retracts to
%$\rho^{-1}(\scrO_{c,N_c})$.
\end{lemma}
\begin{proof}
Let $\Rrr^+ \times U(c,N)$ be the action coming from the cone
contraction given in Lemma~\ref{lemmaUcNccone}.  Since $\rho\colon
\Pee_{\scrF_c(N_c)} \to U(c,N)$ is $T_{N}$-equivariant, this
$\Rrr^+$ action lifts to $\Pee_{\scrF'_c}$. This defines a
deformation of $\Pee_{\scrF'_c}$ into an arbitrarily small
neighborhood of $\rho^{-1}(0_c)$.  Of course, a sufficiently small
neighborhood of $\rho^{-1}(0_c)$ deformation retracts onto
$\rho^{-1}(0_c)$.  Compactness of the fiber follows from the
properness of $\rho$ (see \cite[\S 2.4]{Ful}).
\end{proof}

The refinement $\scrF'_c$ of $\scrF_c$ determines a fan denoted
$\scrF'_c(N_c)$ in $(N_c)_\Rrr$.  As above we have $\rho \colon
\Pee_{\scrF'_c(N_c)} \to V(c)$ with $\Pee_{\scrF'_c(N_c)}$
deformation retracting to the compact variety $\rho^{-1}(0_c)$.
\begin{corollary}\label{prodcor}
$\Pee_{\scrF'_c}$ is a locally trivial fiber bundle over $T_N/T_{N_c}$ with
fiber $\Pee_{\scrF'_c(N_c)}$.  The fibers $\rho^{-1}(\scrO_c) \subset
\Pee_{\scrF'_c}$ and $\rho^{-1}(0_c) \subset \Pee_{\scrF'_c(N_c)}$ are related
by $\rho^{-1}(\scrO_c) \cong \rho^{-1}(0_c) \times T_N/T_{N_c}$.
\end{corollary}
\begin{proof}
Use the fact that $\rho \colon \Pee_{\scrF'_c} \to U(c,N)$ is
$T_N$-equivariant and the fibration structure in
Equation~\ref{UcNcfibration}.
\end{proof}

\subsubsection{Refinements and induced maps: The projective case}

One class of especially nice lattice polytopes are reflexive
polytopes: Given a polytope $\Delta \subset N_\Rrr$, define the
{\em polar polytope}
$$\Delta^\circ = \{ v \in M_\Rrr \, | \, \langle x,v
\rangle \geq -1 , \ \mbox{for all} \ x \in \Delta \} \subset M_\Rrr \ .$$ A
lattice polytope $\Delta$ is said to be {\em reflexive} if $\Delta^\circ$ is a
lattice polytope in $M_\Rrr$.

\begin{definition}
Let $\Delta$ be a reflexive polytope.  Let $\scrT$ be a triangulation which is
a rectilinear subdivision of $\partial \Delta$, with set of vertices equal to
$N \cap \partial \Delta$. The cones of the simplices of $\scrT$ form a
refinement $\scrF(\scrT)$ of $\scrF(\Delta)$.  Such refinements are called {\em
maximal projective subdivisions}.  We denote by $V(\scrT)$ the set of vertices
of the triangulation $\scrT$, i.e., $V(\scrT)=N \cap \pa \Delta$.
\end{definition}

\begin{lemma} \label{LemmaModels}
Let $\Delta$ be a reflexive polytope.  Let $\scrF'$ be a maximal projective
subdivision of $\scrF(\Delta)$.  Let $D_\infty \subset \Pee_{\scrF'}$ be the
Weil divisor $\cup_{v \in V(\scrT)} \overline{\scrO_v}$,
%$\Pee_{\scrF'} \setminus T_N$
called the {\em divisor at infinity}. For any simplex $\delta$ of $\scrT$ of
dimension $\leq 2$, $U( \delta,N)$ is smooth and $D_\infty \cap U(\delta,N)$ is
an anti-canonical divisor in $U(\delta,N)$ and is a divisor with normal
crossings in $U(\delta,N)$.
\end{lemma}
\begin{proof}
By reflexivity of $\Delta$, the vertices $\{v_1, \ldots, v_r\}$ of
any simplex $\delta$ of $\scrF'$ of dimension $\leq 2$ form part
of a basis for $N$ \cite[Corollary A.2.3]{CK}; let $\{v_1,
\ldots, v_r, e_{r+1}, \ldots, e_n\}$ be a basis for $N$.  Use this
basis to split $N = A \oplus B$ with $v_1, \ldots, v_r$ being a
basis for $A$ and $e_{r+1}, \ldots, e_n$ a basis for $B$.  Let
$A_+^*$ denote the ``positive $2^r$-ant'' in the dual space $A^*$
with respect to the dual basis.  Then
$$M \cap \check{c} = A_+^* \times B^*$$
so that
$$U(\delta,N) = \Spec{\Cee}[M \cap \check{c}] = \mathbb{C}_A \times
T_B$$ where $\mathbb{C}_A = \mathbb{C} \otimes_\Zee A$ and
\begin{equation} \label{DinfUeqn}
D_\infty \cap U(\delta,N) = \cup_{i=1}^r H_i \times T_B
\end{equation}
where $H_i \subset \mathbb{C}_A$ is the $i$th coordinate
hyperplane.  Let $t_i : T_N \rightarrow
\mathbb{C}^*$ be the coordinate functions for this basis.  Then
$$\frac{dt_1}{t_1} \wedge \cdots \wedge \frac{dt_n}{t_n}$$
is a holomorphic differential $n$-form on $T_N$ with simple poles
along $D_\infty$ \cite[\S~4.3]{Ful}, proving that $D_\infty$ is
anti-canonical.
\end{proof}

In order to capture all of these properties we make the following definition.
\begin{definition} \label{model}
Let $\Sigma$ be a finite cell complex, not necessarily of
homogeneous dimension.  Let $\scrS$ be a compact complex algebraic
variety. We say that $\Sigma$ is a {\em combinatorial model} for
$\scrS$ if the following hold:
\begin{itemize}
\item For each $i$-cell $c$ of $\Sigma$ there is a locally closed, irreducible
algebraic subset $S_c$ of $\scrS $ and $S _c$ is isomorphic to a
complex torus of dimension $i$. \item $\scrS =\cup_cS _c$. \item
For $c\not= c'$, $S_c\cap S_{c'}=\emptyset$. \item For every cell
$c$ of $\Sigma$ the closure $\overline S_c$ of $S_c$ is isomorphic
to a smooth variety and
$$\overline S_c= \bigcup_{\{c'\prec c\}} S_{c'}.$$
\item For every cell $c$ of $\Sigma$ the subset $\cup_{\{c'\prec
c\left|\right. c'\not=c\}} S_{c'}$ in $S_c$ is a smooth divisor
with normal crossings in $\overline S_c$. It  is an anti-canonical
divisor for $\overline S_c$.
\end{itemize}
\end{definition}

\begin{figure}[h] 
\begin{center}
$\begin{array}{c@{\hspace{.5in}}c}
\includegraphics[width=2in]{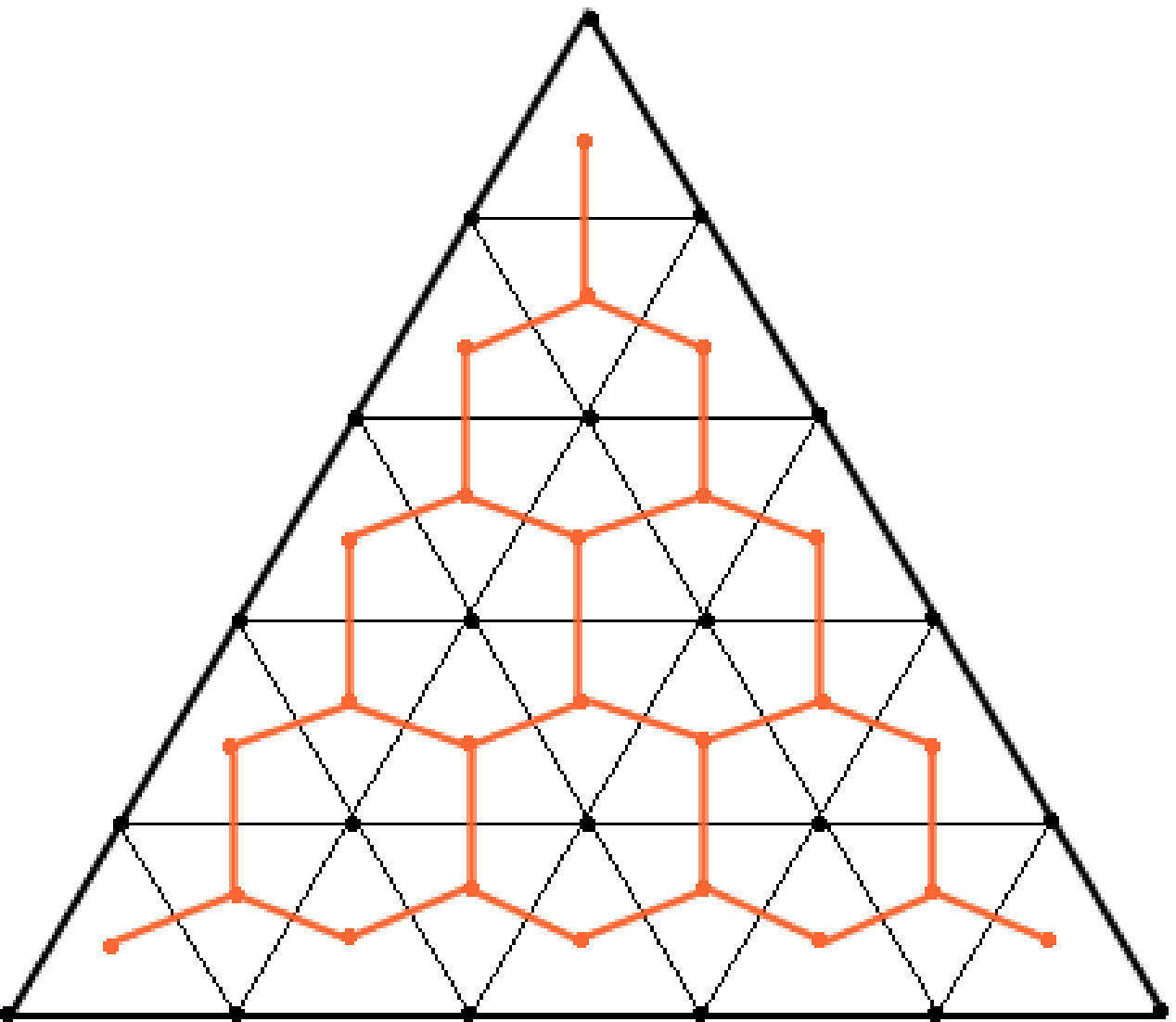} &
\includegraphics[width=2in]{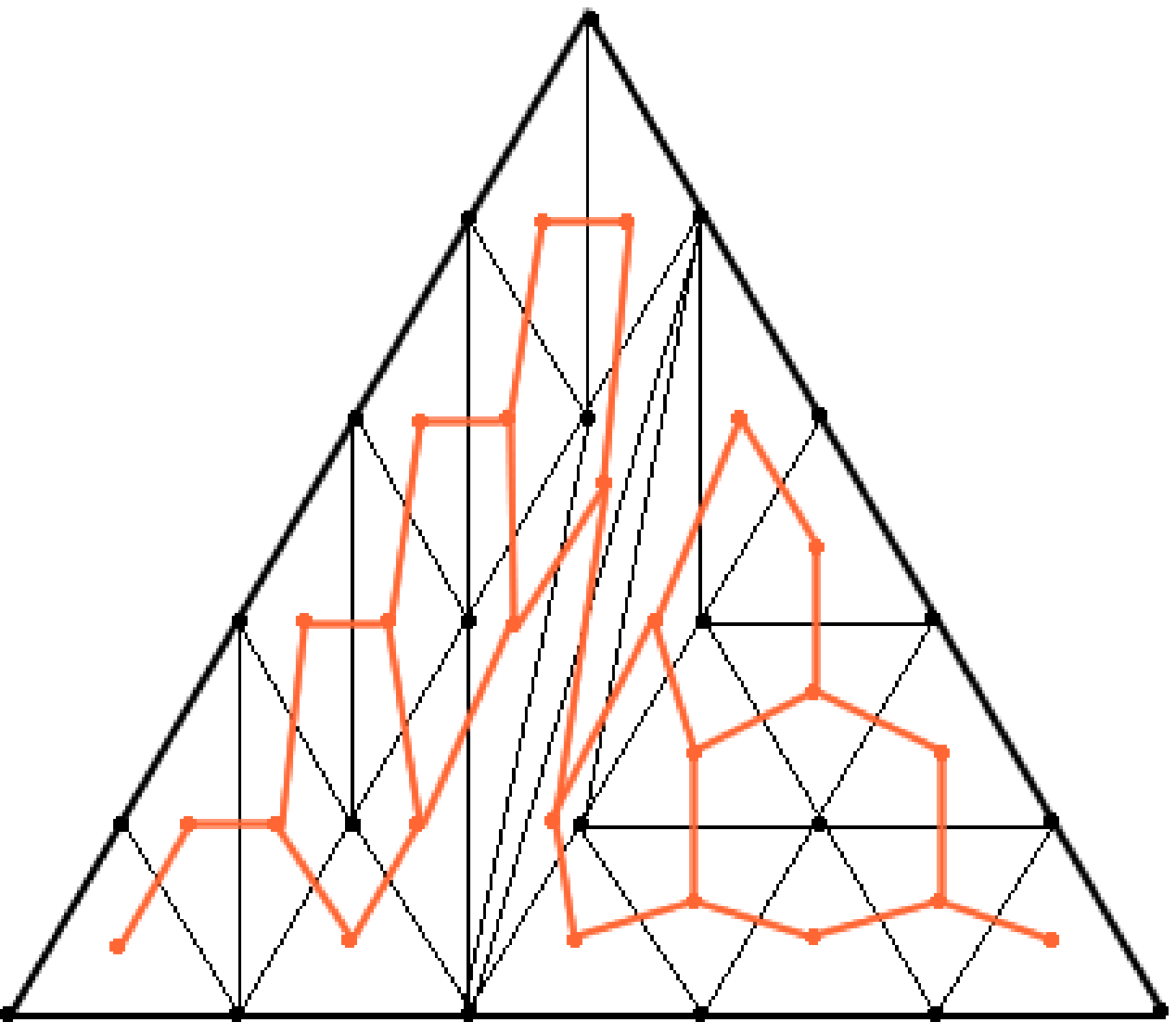}
\end{array}$
\end{center}
\caption{Examples of triangulated two-faces with dual cell
complexes}
\label{figuretriang}
\end{figure}

\begin{prop} \label{propmodels}
Let $\Delta$ be a reflexive polytope, let $\scrF'$ be a maximal
projective subdivision of $\scrF(\Delta)$, and let $f$ be a face
of $\pa \Delta$. Let $\rho\colon \Pee_{\scrF'_f(N)}\to U(f,N)$  be
the map given in Lemma~\ref{LemmaRho} and let $\rho_f$ be the
restriction of $\rho$ to $\Pee_{\scrF'_f(N_f)}\subset
\Pee_{\scrF'_f(N)}$. Denote by $\Sigma(\scrF'_f)$ the cell complex
$\Sigma(\scrF'_{c(f)})$.  Then $\Sigma(\scrF'_f)$ is a
combinatorial model for $\scrS(f)=\rho^{-1}(0_f) \subset
\Pee_{\scrF'_f(N_f)}$. Furthermore, under the map $\rho\colon
\Pee_{\scrF'_f(N)}\to U(f,N)$ the preimage of the closed orbit
$\scrO_f$ is isomorphic to $\scrS(f)\times T_N/T_{N_f}$.
\end{prop}

\begin{proof}
By Lemma \ref{preimage}, $\rho^{-1}(0_f)$ is the union of all the
$T_{N_c}$-orbits in $\scrF'_f(N_f)$ corresponding to simplices of
$\scrF'$ meeting $\stackrel{\circ}{f}$.  The dimension of the orbit
corresponding to a simplex $\delta$ is the codimension of $\delta$
in $f$. Hence there is a bijective correspondence between the
$i$-cells of $\Sigma(\scrF')$ and the $i$-dimensional
$T_{N_f}$-orbits in $\rho^{-1}(0_f)$. By Equation (\ref{face}),
under this correspondence the face relation becomes the adherence
relation.  This establishes the first four conditions in the
definition.  The fifth is immediate from Lemma \ref{LemmaModels}.

The last statement is immediate from Corollary~\ref{prodcor}.
\end{proof}

\subsection{Lefschetz theorems and arithmetic genus of
Calabi-Yau manifolds in toric varieties}

Suppose that $N$ is a lattice of rank $4$ and $\Delta \subset
N_\Rrr$ is a reflexive polytope with $\Pee_\Delta$ as associated
toric variety and $D_\infty$ the divisor at infinity. Since
$\scrO(D_\infty)$ is semi-ample \cite[Lemma~4.1.2]{CK}, the sections of this
bundle define a projective embedding
$$\iota\colon \Pee_\Delta \hookrightarrow \Pee(H^0(\Pee_\Delta,
\scrO(D_\infty))^\vee). $$ Let $\scrT$ be a maximal triangulation
of $\Delta$, i.e., one whose set of vertices is $N \cap \partial
\Delta$. Taking the cones over the simplices of $\scrT$ defines a
maximal projective subdivision $\scrF'=\scrF'(\scrT)$ of
$\scrF(\Delta)$. There are the associated projective variety
$\Pee_{\scrF'}$ and the $T_N$-equivariant morphism
$$\rho \colon \Pee_{\scrF'} \to \Pee_\Delta.$$
Recall that for each face $\gamma$ of $\partial \Delta$,
$\scrF'_\gamma(N_\gamma)$ is the sub-fan of $\scrF'$ consisting of all cones in
$\scrF'$ contained in $c(\gamma)$. According to \cite[Lemma 4.1.2 and Corollary
A.2.3]{CK} for every face $\gamma$ of $\partial \Delta$ of positive
codimension, $\Pee_{\scrF'_\gamma(N_\gamma)}$ is smooth.

Let $Y \subset \Pee_\Delta$ be the intersection of $\Pee_\Delta$
with a generic hyperplane in $\Pee(H^0(\Pee_\Delta,
\scrO(D_\infty))^\vee)$. Then $Y$ is a generic anticanonical
divisor in $\Pee_\Delta$ transverse to all $T_N$-orbits. We set
$\widetilde{Y} = \rho^{-1}(Y) \subset \Pee_{\scrF'}$. It is also
anti-canonical and, as was the case with
$\Pee_{\scrF'_\gamma(N_\gamma)}$ above, $\tY$ is a smooth
Calabi-Yau threefold. For any lattice point $\ell\in N\cap
\pa\Delta$, there is a toric divisor $\widetilde \scrO_\ell\subset
\Pee_{\scrF'}$. The only singularities of its closure in
$\Pee_{\scrF'}$ occur at fixed points of the $T_N$-action. Thus,
the intersection of $\tY$ with the closure of
$\widetilde\scrO_\ell$ is a smooth divisor in $\tY$. We denote it
by $R(\ell)$. Its image in $Y$ is the intersection of $Y$ with the
closure of the orbit corresponding to the unique open face of $\pa
\Delta$ containing $\ell$.

Let us turn to the analogue for $\Delta$ of dimension $n > 4$ and
complete intersections. First we need to recall the notion of
Batyrev and Borisov of NEF partitions defining complete
intersections. Let $\Delta\subset N_\Rrr$ be a reflexive polytope.
Let $V(\Delta)=V_1\coprod \cdots \coprod V_k$ be a partition of
the vertices of $\Delta$. According to \cite[\S~4]{BB} this is a {\em
NEF partition} if there are functions $\varphi_i\colon N_\Rrr\to
\Rrr,\ 1\le i\le k$, with the following properties:
\begin{enumerate}
\item $\varphi_i$ is piecewise linear and linear on the cone on
each face of $\Delta$. \item $$\varphi_i(v)=\begin{cases} 0 & \ \
{\rm if\ } v\not\in V_i \\  -1 & \ \ {\rm if\ } v\in
V_i\end{cases}.$$ \item $\varphi_i$ is a concave function; i.e.,
$\varphi_i(ta+(1-t)b)\ge t\varphi_i(a)+(1-t)\varphi_i(b)$. \item
For each $i\le k$ and each maximal dimensional face $f$ of
$\partial \Delta$ there is an $m_{i,f}\in M$ such that
$$\varphi_i|_{c(f)}=\langle m_{i,f},\cdot\rangle.$$
\end{enumerate}
Set $\nabla_i^\circ\subset M_\Rrr$ equal to the convex hull of the
$\{m_{i,f}\}$ as $f$ ranges over the maximal dimensional faces of
$\partial \Delta$. Then
$$\varphi_i(x)=\min_{y\in
\nabla_i^\circ}\langle y,x\rangle.$$

Let $D_i=\sum_{v\in V_i}\overline{\scrO}_{v}$. Clearly,
$D_\infty=D_1+\cdots +D_k$ is anti-canonical. The divisor $D_i$ is
ample if and only if the function $\varphi_i$ is strictly convex
on $\partial \Delta$, see \cite[p.~70]{Ful}. Now let us suppose
that we have a NEF partition $V(\Delta) = V_1\coprod \cdots
\coprod V_k$ with the property that each of the corresponding
divisors $D_i$ is ample.  This condition implies that
$\nabla_i^\circ$ is combinatorially dual to $\Delta$.

We need the following elementary lemma from \cite[Proposition 6.3]{BB}.

\begin{lemma}\label{hypersurface}
Suppose that $\gamma$ is a face of $\partial \Delta$ and suppose that the
relative interior of $\gamma$ contains a lattice point. Then there is an $i\le
k$ such that all vertices of $\gamma$ belong to the same $V_i$. We denote by
$i(\gamma)$ this index.
\end{lemma}

Suppose that the relative interior of $\gamma$ contains a lattice
point.  Then all the vertices of $\gamma$ are contained in
$V_{i(\gamma)}$, and hence for every $j\not= i(\gamma)$ we define
the dual of $\gamma$ in $\nabla_j^\circ$, denoted $\gamma^\vee_j$,
to be the face of $\nabla_j^\circ$ that evaluates $0$ on $\gamma$.
Analogously, we define $\gamma^\vee_{i(\gamma)}$ to be the face of
$\nabla_{i(\gamma)}^\circ$ that evaluates $-1$ on $\gamma$. Since
$D_i$ is ample, $\gamma_i^\vee$ is the face of $\nabla_i^\circ$
dual to the face $\gamma$ of $\pa \Delta$ under the combinatorial
duality. In particular, for each face $\gamma$ of $\pa\Delta$ the
dimension of $\gamma_i^\vee$ is $n-1-{\rm dim}(\gamma)$.

\subsubsection{Lefschetz theorems} \label{lefthmsec}

We need a generalization of the classical Lefschetz theorem, a
generalization established by Bernstein, Danilov and Khovanskii
(\cite[Theorem 6.4]{DK}). Consider the generic affine complete
intersection $Y_0 \subset {\mathbb{T}}^n = (\Cee^*)^n$ of
dimension $r=n-k$.
\begin{theorem} \label{LefThm}
The compactly supported cohomology $H^i_c(Y_0) = 0$ for $i < r$. Furthermore,
if all the Newton polyhedra $\Delta_1, \ldots, \Delta_k$ have dimension $n$,
then the Gysin homomorphism $H^i_c(Y_0) \to H^{i+2k}_c({\mathbb{T}}^n)$ is an
isomorphism for $i > r$ and surjective for $i=r$.
\end{theorem}

The first statement is classical, given that $Y_0$ is a smooth affine variety
of dimension $r$.

Applying Poincar\'{e} duality for noncompact manifolds yields:
\begin{corollary} \label{LefCor}
Under the hypotheses of Theorem~\ref{LefThm}, the map in homology
induced by the inclusion,
$$H_i(Y_0) \to H_i({\mathbb{T}}^n)$$
is an isomorphism for $i$ less than the complex codimension $k$ of
$Y_0$ and surjective for $i=k$.
\end{corollary}

It will be important to have these results not only for the open torus-orbit in
$P_\Delta$ but also for orbits of smaller dimension. In the hypersurface case,
let $\omega$ be a generic section of $\scrO(D_\infty)$. Consider a face $\gamma
$ of $\Delta$. The restriction of $\omega$ to the orbit $\scrO_\gamma$ has
convex hull an affine translate of $\gamma^\vee\subset \Delta^\circ$ and hence
the dimension of this convex hull is equal to the dimension of the orbit
$\scrO_\gamma$. Thus, the Lefschetz theorem applies to $Y\cap \scrO_\gamma\to
\scrO_\gamma$. But notice that  the dimension of $Y\cap \scrO_\gamma$ is
smaller than $r$, and the comparison of homologies only goes up to the middle
(real) dimension of this variety.

In the complete intersection case, the necessary dimension hypothesis does not
automatically restrict to the faces. For this we need the $D_i$ to be ample.

\begin{lemma}
Let $\Delta\subset N_\Rrr$ be a reflexive polytope and let
$V(\Delta)=V_1\coprod\cdots\coprod V_k$ be a NEF partition. Suppose that each
of the divisors $D_i$ associated to the $V_i$ are ample. Let $\omega_i$ be a
generic section of $D_i$. Then for each face $\gamma$ of $\pa\Delta$, the
convex hull of the support of $\omega_i|_\gamma$ has dimension equal to the
dimension of $\scrO_\gamma$.
\end{lemma}

\begin{proof}
The convex hull of the support of $\omega_i|_\gamma$ is an affine translate of
$\gamma_i^\vee\subset \nabla_i^\circ$. Since $D_i$ are ample, the
$\nabla_i^\circ$ are combinatorially dual to $\Delta$, and hence the faces
$\gamma_i^\vee$ have  dimension $n-1-{\rm dim}(\gamma)$, which is exactly the
dimension of $\scrO_\gamma$.
\end{proof}

Now applying Corollary~\ref{LefCor} gives the following:

\begin{prop}\label{lefcor}
Let $\Delta\subset N_\Rrr$ be a reflexive polytope. Let $Y\subset P_\Delta$ be
an $r$-dimensional variety, either a generic anti-canonical hypersurface or a
generic complete intersection associated to a NEF partition
$D_\infty=D_1+\cdots+D_k$, with the $D_i$ being ample. Then for each face
$\gamma$ of $\Delta$ the inclusion $Y\cap\scrO_\gamma\to\scrO_\gamma$ induces
an isomorphism on homology below  dimension $(r-{\rm codim}(\gamma))$ and
induces a surjection on homology in that dimension.
\end{prop}

\subsubsection{Arithmetic genus}

We begin by recalling a result of Khovanskii (\cite[Theorem
1]{Kho}) in the case of hypersurfaces in a torus.
 \begin{prop}
 Let $\Cee^*\otimes N$ be a torus, and let $\omega\in \Cee[\mathbb T]$ be a
 non-degenerate regular function on $\Cee^*\otimes N$, non-degenerate in the sense that
 $Y=\omega^{-1}(0)$ is a smooth hypersurface.
 Let $\Delta'\subset M_\Ar$ be the convex hull of the support of $\omega$.
 Then the arithmetic genus $\chi(Y)$ is given by
$$\chi(Y) = 1 - B(\Delta')$$
where
$$B(\Delta') = (-1)^{\dim(\Delta')} \#(\Delta' \cap M) \ .$$
\end{prop}

We apply this to our context. First, let us also suppose that $N$
is of rank $4$ and that $\Delta\subset N\otimes \Ar$ is reflexive.
Let $\gamma$ be a face of $\partial \Delta$ and denote by
$\gamma^\vee$ the dual face in $\Delta^\circ$ consisting of all
$x\in \Delta^\circ$ with the property that $\langle
x,\gamma\rangle =-1$. Under the mapping $\Pee_\Delta\to
\Pee(H^0(\Pee_\Delta,\scrO(D_\infty))^\vee)=\Pee(\Cee^{M\cap
\Delta^\circ})$ the image of the orbit $\scrO_\gamma$ is contained
in the subprojective space determined by the vanishing of all
$m\in \Delta^\circ \setminus \gamma^\vee$. Thus, the supporting
polytope for the restriction of a generic section
$$\omega = \sum_{m \in M \cap \Delta^\circ} a_m \chi^m$$
of $\scrO(D_\infty)$ to $\scrO_\gamma$ is the convex hull of $M
\cap \gamma^\vee$.  Since the vertices of $\gamma^\vee$ are
contained in $M$, the supporting polytope is $\gamma^\vee$.
Applying the Khovanskii result cited above to $\scrO_\gamma$ gives
the following:

\begin{corollary}\label{dim1&2hyper}
Let $N$ be a lattice of rank 4, $\Delta \subset N_\Rrr$ a
reflexive polytope, $\omega$ a generic section of
$\scrO(D_\infty)$, and $Y\subset \Pee_\Delta$ the zero locus of
this section. Then for any edge $e$ of $\Delta$ the intersection
$Y\cap \scrO_e$ is a (non-compact) riemann surface of genus
$\ell^*(e^\vee)$. For any two-face $f$ of $\Delta$ the
intersection $Y\cap \scrO_f$ consists of $1+\ell^*(f^\vee)$
points.
\end{corollary}

Let's now generalize this to the complete intersection case.

Theorem 1 of \cite{Kho} also implies the following.

\begin{prop}
The arithmetic genus $\chi(Y)$ of the variety $Y$ defined in
$(\Cee^*)^n$ by a nondegenerate system of equations $\omega_1 =
\ldots = \omega_k = 0$ with Newton polyhedra $\Delta'_1, \ldots,
\Delta'_k$ is given by
$$\chi(Y) = 1 - \sum_i B(\Delta'_i) + \sum_{i<j} B(\Delta'_i +
\Delta'_j) - \ldots + (-1)^k B(\Delta'_1 + \ldots + \Delta'_k) \
,$$ where
$$B(\Delta') = (-1)^{\dim(\Delta')} \#(\Delta' \cap M) \ .$$
\end{prop}

Applying this result to the various faces gives:

\begin{corollary}\label{dim1&2ci}
Let $\Delta \subset N_\Rrr$ be a reflexive polytope, where $\dim
N_\Rrr = n = k+3$.  Let
$$V(\Delta)=\coprod_{i\in I} V_i \ ,$$
where $I = \{1,\ldots,k\}$,
 be
a NEF partition with the associated divisors $D_i$ being ample.
Let $Y \subset \Pee_\Delta$ be a generic complete intersection of
sections of the $\scrO(D_i)$.  Then: \begin{enumerate} \item  For
any 2-face $f$ of $\pa \Delta$ whose relative interior contains a
lattice point,
$$\# (Y \cap \scrO_f) =1+\sum_{J \subset I;\ J\not=\emptyset} (-1)^{3 -
|J|} \ell^*( \sum_{j \in J} f_j^\vee ).$$  \item For an edge $e$
of $\pa \Delta$ containing an interior lattice point, the curve $Y
\cap \scrO_e$ has
$$\chi(Y \cap
\scrO_e) = 1-\sum_{J \subset I;\ J\not=\emptyset} (-1)^{r - |J|}
\ell^*(\sum_{j \in J} e^\vee_j).$$
\end{enumerate}
\end{corollary}

Notice that when $k=1$, this formula specializes to the one given in
Corollary~\ref{dim1&2hyper}.

\subsection{First computations of homology groups of Calabi-Yau
threefolds} \label{subsectfirst}

Throughout this section we suppose that $\Delta$ is a reflexive
polytope and that either the dimension of $\Delta$ is $4$ and $Y$
is the vanishing locus of a generic section od $\scrO(D_\infty)$,
or that the dimension of $\Delta$ is $n=k+3$ and we have a NEF
partition $V(\Delta)=V_1\coprod\cdots\coprod V_k$ with the
corresponding divisors $D_1,\ldots,D_k$ being ample with $Y$ being
the complete intersection of the zero loci of generic sections of
the $D_i$.  In either case $Y$ is a (possibly singular) Calabi-Yau
three-fold. Since $\tY=\rho^{-1}(Y) \subset\Pee_\scrT$ is smooth,
it is easy to compute $H_1(\tY)$ and $H_2(\tY)$ given the
Lefschetz theory described above.

Let $\tY_0$ be the intersection of $\tY$ with the open $T_N$-orbit in
$\Pee_\scrT$.
 By the
Lefschetz theorem (Section~\ref{lefthmsec}) we have
\begin{eqnarray*}
H_0(\tY_0)=H_0(Y_0) &  =  & \Zee \\
H_1(\tY_0)=H_1(Y_0) & = &  N \\
H_2(\tY_0)=H_2(Y_0)  & =  & \wedge^2 N. \end{eqnarray*}

We set $\tW\subset \tY$ equal to the union of $Y_0$ and all the intersections
of $\tY$ with the codimension-one torus orbits in $\Pee_\scrT$. Of course,
these orbits are indexed by the vertices of $\scrT$, i.e., by $N\cap \partial
\Delta$. Thus, the pair $(\tW,\tY_0)$ is excisively equivalent to
$$\coprod_{\ell\in \partial\Delta\cap N}(\tY \cap \widetilde\scrO_\ell) \times (D^2, S^1) .$$
Let's use the Lefschetz theory to analyze $\tY\cap \widetilde\scrO_\ell$.

\begin{lemma}\label{2.22}
\begin{enumerate}
\item If $\ell\in V(\partial \Delta)$, then $\tY\cap \widetilde
\scrO_\ell=Y\cap\scrO_\ell$. Furthermore, $H_0(\tY\cap\widetilde
\scrO_\ell)=\Zee$ and the inclusion $\tY\cap
\widetilde\scrO_\ell\to \widetilde\scrO_\ell$ induces an
isomorphism $H_1(\tY\cap \widetilde \scrO_\ell)\to
H_1(\widetilde\scrO_\ell)$. \item If $\ell$ is contained in the
interior of an edge $e$ of $\Delta$, then $\tY\cap \widetilde
\scrO_\ell$ is a torus-fibration over $Y\cap \scrO_e$ with
(one-dimensional) fiber
$$\frac{N\cap\Rspan(e)}{N\cap \Rspan(\ell)}\otimes \Cee^*.$$
Furthermore, $\tY\cap\widetilde \scrO_\ell$ is connected and the
inclusion $\tY\cap \widetilde\scrO_\ell\to \widetilde \scrO_\ell$
induces a surjection on $H_1$. \item If $\ell$ is contained in the
interior of a two-face $f$ of $\Delta$, then $\tY\cap\widetilde
\scrO_\ell$ is a torus fibration over $Y\cap \scrO_f$, which is a
finite, non-empty set of points. The fiber is identified with the
complex two-torus
$$\frac{N\cap\Rspan(f)}{N\cap \Rspan(\ell)}\otimes \Cee^*.$$
In particular, each component of $\tY\cap\widetilde\scrO_\ell$ has
first homology isomorphic to
$$\frac{N\cap\Rspan(f)}{N\cap \Rspan(\ell)}$$
inside $H_1(\widetilde\scrO_\ell)=N/\langle \ell\rangle$. \item If
$\ell$ is contained in the interior of a three-face of $\partial
\Delta$, then $\tY\cap \widetilde \scrO_\ell=\emptyset$.
\end{enumerate}
\end{lemma}

\begin{proof}
If $\ell$ is contained in the interior of a face $\gamma$ of $\Delta$, then the
restriction of the map $\rho$ to a map $\widetilde\scrO_\ell\to\scrO_\gamma$ is
the natural projection
$$\frac{N}{\Zee\langle \ell\rangle}\otimes \Cee^*\to
\frac{N}{\Rspan(\gamma)}\otimes \Cee^*.$$
Since $Y \subset X$ is generic, $\tY \cap \scrO_\gamma$ is a fibration with
fiber $$\frac{\Rspan(\gamma)}{\Zee\langle \ell\rangle} \otimes \Cee^* .$$

When $\ell$ is a vertex of $\pa\Delta$, $\tY \cap \widetilde\scrO_\ell = Y \cap \scrO_\ell$
and the homology statements in part (1) are verified in
Proposition~\ref{lefcor}.

Let us suppose that $\ell$ is contained in the interior of an edge $e$ of
$\Delta$. Then by Proposition~\ref{lefcor} we have $Y\cap \scrO_e$ is connected
and $H_1(Y\cap\scrO_e)\to H_1(\scrO_e)$ is surjective. Since
$\tY\cap\widetilde\scrO_\ell\to \widetilde\scrO_\ell$ is a $\Cee^*$-bundle over
the inclusion $Y\cap\scrO_e\to\scrO_e$, the result follows in this case.

Now suppose that $\ell$ is contained in the interior of a two-face $f$. By
Proposition~\ref{lefcor}, $Y\cap \scrO_f$ is a finite, non-empty, set of
points. The result is then clear in this case.

Lastly, if $\ell$ is contained in the interior of a three-face
$g$, then the generic $Y$ in the given linear system misses the
point $\scrO_g$ in $P_\Delta$. Hence, the preimage $\tY$ misses
$\widetilde\scrO_\ell$.
\end{proof}

 Using the long exact sequence of the pair $(\tW,\tY_0)$  gives
$$ \hspace{-.3in} \bigoplus_{\ell \in \pa \Delta^{(2)} \cap N} H_1(\widetilde\scrO_\ell \cap
\tY) \stackrel{\oplus_\ell (\cdot \wedge \ell)}{\longrightarrow}
\wedge^2 N \to H_2(\tW) \to \bigoplus_{\ell \in \pa \Delta^{(2)}
\cap N} H_0(\widetilde\scrO_\ell \cap \tY)
\stackrel{\bigoplus_\ell\cdot \ell}{\to} N \to H_1(\tW) \to 0.$$
Here, the notation $\pa\Delta^{(2)}$ refers to the two-skeleton of
$\pa\Delta$. The first map is the direct sum over $\ell$ of the
compositions
$$H_1(\tY \cap \widetilde\scrO_\ell)\to H_1(\widetilde\scrO_\ell)=N/\langle
\ell\rangle\stackrel{\cdot\wedge\ell}{\to}\wedge^2N.$$ Using this
we shall show the following result.

\begin{corollary}\label{2.23}
\begin{enumerate}
\item
$H_1(\tY)=N/{\rm Span}(N\cap \partial \Delta^{(2)})$. \item We have an exact
sequence:
$$0 \to \Tor H_2(\tY) \to H_2(\tY) \to \Ker (\bigoplus_\ell(\cdot \ell)) \to 0.$$
\item $\Ker (\bigoplus_\ell(\cdot \ell))$ is free
abelian. In the hypersurface case it is of rank
$$\rk \Ker (\bigoplus_\ell(\cdot \ell)) = \# V + \sum_e \ell^*(e) + \sum_f
\ell^*(f) (1 + \ell^*(f^\vee)) - 4 .$$ In the complete intersection
case (of ample divisors) its rank is
$$\hspace{-.05in} \rk \Ker(\bigoplus_\ell(\cdot \ell)) = \# V + \sum_e \ell^*(e) + \sum_f
\ell^*(f)\left(1+\sum_{J \subset I;\ J\not=\emptyset} (-1)^{3 -
|J|} \ell^*( \sum_{j \in J} f_j^\vee)\right)-n .$$ \item
$$\hspace{-.25in} \Tor H_2(\tY)  = \wedge^2 N / \Image \left( \bigoplus_{\ell \in
\pa \Delta^{(1)} \cap N} \frac{N}{N\cap \Rspan(\ell)} \wedge \ell
\oplus \bigoplus_f\left(\bigoplus_{\ell\in
\stackrel{\circ}{f}}\frac{N\cap
\Rspan(f)}{N\cap\Rspan(\ell)}\wedge \ell\right)\right).$$
\end{enumerate}
\end{corollary}

\begin{proof}
Since the complement of $\tW$ in $\tY$ is a union of orbits of complex
codimension at least two, and since $\tY$ is smooth, it follows by general
position that $H_i(\tW)\to H_i(\tY)$ is an isomorphism for $i\le 2$. Hence, we
work with $H_i(\tW)$ for $i=1,2$.

The first statement is immediate from the long exact sequence of the pair
$(\tW,\tY_0)$ since $\widetilde\scrO_\ell\cap \tY$ is non-empty for every
$\ell\in
\partial\Delta^{(2)}$. It is easy to see that the first map in this exact
sequence has torsion cokernel. Since the fourth term is free
abelian, the second statement follows. It also follows that
$\Ker(\bigoplus_\ell(\cdot\ell))$ is free abelian. Let us compute
its rank, which is the rank of fourth term minus $n$, the rank of
$N$. By the Lefschetz theorem $\widetilde\scrO_\ell\cap \tY$ is
connected if $\ell$ is contained in the one-skeleton of $\partial
\Delta$. By Corollary~\ref{dim1&2hyper} in the hypersurface case
$\widetilde\scrO_\ell\cap \tY$ has $(1+\ell^*(f^\vee))$ components
if $\ell$ is contained in the relative interior of two-face $f$. By
Corollary~\ref{dim1&2ci} in the complete intersection case if $\ell$
is contained in the relative interior of a two-face then
$\widetilde\scrO_\ell\cap \tY)$ has
$$1+\sum_{J \subset I;\ J\not=\emptyset} (-1)^{3 - |J|} \ell^*(
\sum_{j \in J} f_j^\vee)$$  components. In both cases the formula
for the rank is immediate. The fourth item is clear from our
description of $H_1(\widetilde\scrO_\ell\cap\tY)$.
\end{proof}

Notice that parts (2) and (3) of this result establish part (1) of
Theorem~\ref{thmH2intro} and gives the more general formula for the
rank of $H_2(\tY)$ in the complete intersection case.

\subsubsection{Representation by divisors}

We now establish the results in part (2) of Theorem~\ref{thmH2intro}.

\begin{prop}
For every $\alpha \in H^2(\tY)$, there exists a integral linear combination of
divisors  $D=\sum_jn_jD_j$ such that $\alpha$ is Poincar\'e dual to $D$. In
fact, each irreducible component $D_j$ in this sum can be taken to be an
irreducible component of the intersection of $\tY$ with the closure of a
torus-orbit of codimension one in $\tX$.
\end{prop}

\begin{proof}
We have the exact sequence
$$\bigoplus_\ell H^0(\tY\cap \widetilde\scrO_\ell)\to H^2(\tY)\to H^2(\tY_0).$$
The image of the generator of one of the $\Zee$-summands in
$H^0(\tY\cap\widetilde \scrO_\ell)$ maps to the Poincar\'e dual of
the corresponding divisor in $H^2(\tY)$. Thus, the proposition
follows immediately if we can show that the map $H^2(\tY)\to
H^2(\tY_0)$ is trivial. Since $H^2(\tY_0) = H^2(T_N)$ is torsion
free, it suffices to show that the algebraically dual map
$$H_2(\tY_0)\to H_2(\tY)$$ has torsion image.
Of course, the inclusion $\tW\subset \tY$ induces an isomorphism
$H_2(\tW)\to H_2(\tY)$, and hence we need only see that
$$\bigoplus_\ell H_1(\widetilde \scrO_\ell\cap\tY)\to \wedge^2N=H_2(\tY_0)$$
 $\wedge^2N = H_2(\tY_0)\to H_2(\tY)$ has torsion image. This we already
 observed in the proof of Corollary~\ref{2.23}

 This completes the
proof.
\end{proof}

\begin{corollary}
The Hodge structure on $H^2(\tY)$ is of type $(1,1)$.
\end{corollary}

\begin{corollary}
The mixed Hodge structure of $H^2(Y)$ is pure of weight 2 and Hodge type
$(1,1)$.
\end{corollary}
\begin{proof}
$H^2(Y) \to H^2(\tY)$ is injective.
\end{proof}

Unfortunately, this direct approach we employed here is not as
useful for computing $H_3(Y)$, since we need to understand the
role of the codimension-two orbits. It is also not as useful for
computing the map $H_2(\tY)\to H_2(Y)$ since $Y$ is singular. We
find it convenient to organize the computation differently in
order to address these issues.

\section{Integral homology of Calabi-Yau threefolds in toric varieties}

In this section we compute the maps $H_*(\tY;\Zee)\to H_*(Y;\Zee)$
for $*\le 3$. As we indicated at the end of the last section the
direct approach inductively studying the intersection of $\tY$ with
the orbits in $\Pee_\scrT$ of higher and higher codimension is not
the best way to organize the argument. Rather one inductively
considers the preimage of the intersection of $Y$ with the
torus-orbits in $X$ of higher codimension. The reason is that the
preimage of these intersection is a union of intersections of $\tY$
with orbits in $\Pee_\scrT$ of various dimension, but with the help
of the combinatorial models for these preimages we are able to say a
lot about their homology that is not apparent studying one orbit at
a time in $\tY$.

\subsection{The general set-up}

Let $\Delta$ be a reflexive polytope of dimension $n$. Let $X =
\Pee_\Delta$. Let $X_i$ denote  the union of the $T_N$-orbits of
dimension $\geq n-i$.  We have the chain of inclusions
$$X_0 \subset X_1 \subset X_2 \subset X_3 \subset X_4 \subset \cdots \subset X_n=X .$$
 In
particular $X_0$ is the open torus orbit and hence is isomorphic to
$T_N$. The difference $X_i \setminus X_{i-1}$ is the union of orbits
of complementary dimension, i.e.,
$$X_i \setminus X_{i-1} = \coprod_{\codim \scrO = i} \scrO .$$
Denoting $\Pee_{\scrF'}$ by $\tX$ and by $\rho : \tX \to X$ the
induced resolution, we define $\tX_i = \rho^{-1}(X_i)$ so that we
have
$$\tX_0 \subset \tX_1 \subset \tX_2 \subset \tX_3 \subset \tX_4 \subset \cdots
\subset \tX_n= \tX .$$

As always $Y\subset X$ is either a generic section of
$\scrO(D_\infty)$ if $\Delta$ has dimension $4$ or $Y$ is a
generic complete intersection arising from a NEF partition
$V(\Delta)=V_1\coprod\cdots\coprod V_k$ if $\Delta$ has dimension
$3+k$. In the later case we assume each of the associated divisors
$D_i$ are ample. In both cases we define the nested sequence of
open subvarieties of $Y$
$$Y_0 \subset Y_1 \subset Y_2 \subset Y_3 = Y$$
where
$$Y_i \, = Y\cap X_i \ .$$
Let $Z_i$ be the intersection of $Y$ with all torus orbits of
codimension $i$. Then $$Z_i = Y_i \setminus Y_{i-1}$$ is a closed
subvariety of $Y_i$.  It is a disjoint union $$Z_i = \coprod_{\{f
| \dim(f) = i-1\}} Z(f) \ ,$$ where $Z(f)$ is $Y \cap \scrO_f$. We
let $\nu_i\subset Y_i$ be a (closed) regular neighborhood of
$Z_i$. We choose $\nu_i$ so that the projection $\pi\colon
\nu_i\to Z_i$ is a proper, locally trivial fibration with compact
fibers.  The intersection of $\nu_i$ with any fiber is identified
with a regular neighborhood of the cone point $0_c$ in $V(c)$.
Then $\nu_i^* = \nu_i \cap Y_{i-1}$ is the complement of $0_c$ in
this neighborhood, and $\partial\nu_i$ is a deformation retract of
$\nu_i^*$. Of course, we have
$$Y_i=\nu_i\cup_{\nu_i^*}Y_{i-1}.$$
We have
$$\nu_i=\coprod_{\{f\left|\right. {\rm
dim}(f)=i-1\}}\nu(f) \ ,$$ where $\nu(f)$ is the component of
$\nu_i$ containing $\scrO_f\cap Y$.  We see that the inclusion
induces an identification
$$H_*(\nu(f), \pa \nu(f)) = H_*(\nu(f), \nu^*(f)).$$
Furthermore, each of these pairs is a locally trivial relative
fibration over $Z(f)$ with fiber homotopy equivalent  to the pair
$(U(f,N_f), U^*(f,N_f))$. The $T_N$-action produces a trivialization
of this relative fiber bundle. The relative homology is then
\begin{eqnarray*}
H_*(\nu_i, \pa \nu_i) & \simeq & \bigoplus_{{\rm dim}f=i} \bigoplus_{a+b=*}
H_a(Z(f)) \otimes H_b(V(f),V^*(f)) \\
& = & \bigoplus_{{\rm dim}f=i} \bigoplus_{a+b=*} H_a(Z(f)) \otimes
H_{b-1}(V^*(f)).
\end{eqnarray*}

We define $\widetilde{Y}_i=\rho^{-1}(Y_i)$ for $i=1,\cdots,3$, so
that $\tY_i-\tY\cap\tX_i$. We obtain
$$\widetilde{Y}_0\subset \widetilde{Y}_1\subset \widetilde{Y}_2\subset\widetilde{Y}_3=\tY.$$
Let $\widetilde \nu_i=\rho^{-1}(\nu_i)$,  $\widetilde
\nu_i^*=\rho^{-1}(\nu_i^*)$. We denote by $\widetilde\nu(f)$ and
$\widetilde\nu^*(f)$ the preimages of $\nu(f)$ and $\nu^*(f)$
respectively. We also  denote by $\rho_i\colon \widetilde Y_i \to
Y_i$ the map induced by $\rho$. The pair
$(\widetilde{\nu}(f),\widetilde{\nu}^*(f)) \to Z(f)$ is a relative
fibration with fiber homotopy equivalent  to the pair
$$\left(\Pee_{\scrF'_f(N_f)},\Pee_{\scrF'_f(N_f)}\setminus\scrS(f)\right)$$
where $\scrS(f)$ is as given in Proposition~\ref{propmodels}.
Again by Corollary~\ref{prodcor} the $T_N$-action produces a
trivialization of this bundle.

\subsection{The comparison of $\tY_1$ and $Y_1$}
 We know that $\rho_1 \colon \widetilde{Y}_1 \to Y_1$ is
an isomorphism.

\subsection{The comparison of $\tY_2$ and $Y_2$}

To make this comparison, we need to understand the nature of
$\widetilde \nu_2$ or equivalently, $\Pee_{\scrF'_e(N_e)}$. For any
edge $e$ of $\Delta$,
$$\rho^{-1}(\scrO_e) \simeq \scrO_e \times \scrS(e) \ .$$
Consider the one-cell complex dual to the restriction of $\scrT$ to
${\stackrel{\circ}{e}}$: it has a one-cell for each vertex in
$\stackrel{\circ}{e}$ and a zero-cell for each edge of $\scrT|_e$.
By Proposition~\ref{propmodels} this dual cell complex is a
combinatorial model for $\scrS(e)$ in the sense of
Definition~\ref{model}.
 Thus, $\scrS(e)$ is a
chain of $\Pee^1$'s, one irreducible component for each one-cell.
Two irreducible components intersect if and only if the
corresponding one-cells share a vertex, in which case the
irreducible components meet in the single point which is the
$0$-dimensional orbit corresponding to this vertex. The smooth
affine complex surface $\Pee_{\scrF'_e(N_e)}$ deformation retracts
onto $\scrS(e)$, and hence $H_2(\Pee_{\scrF'_e(N_e)})$ is
identified with $\scrA_{\ell^*(e)}$, the root lattice of the Lie
algebra $A_{\ell^*(e)}$. The intersection pairing on $H_2$ of this
surface agrees with the usual symmetric pairing on this root
lattice given by the Cartan matrix. In particular, the pairing is
non-degenerate and its adjoint has cyclic quotient of order
$k(e)=\ell^*(e)+1$. This means that
$H_2(\Pee_{\scrF'_e(N_e)},\partial \Pee_{\scrF'_e(N_e)})$ is
identified with the dual lattice $\scrA^*_{\ell^*(e)}$ and the
natural map
$$H_2(\Pee_{\scrF'_e(N_e)})\to  H_2(\Pee_{\scrF'_e(N_e)},\partial
\Pee_{\scrF'_e(N_e)})$$ is injective with cokernel, denoted $Q_e$, a
cyclic group of order $k(e)$, see Example~\ref{eg:edge}. Of course,
$$(\widetilde{\nu}_2,\widetilde{\nu}^*_2) \simeq \coprod_{e\in E} Z(e)
\times
\left(\Pee_{\scrF'_e(N_e)},\Pee_{\scrF'_e(N_e)}\setminus\scrS(e)\right)
\ .$$

Define
$$\begin{array}{rclcrcl}
\nu'_2 &=& \overline{\nu_2 \setminus (\nu_2 \cap \nu_3)} &,&
\pa_{\rm hor} \nu'_2 &=& \nu_2 \cap \pa
\nu_3 \\
\widetilde{\nu}'_2 &=& \rho^{-1} (\nu'_2) &,& \pa_{\rm hor}
\widetilde{\nu}'_2 &=& \rho^{-1}(\pa_{\rm hor} \nu'_2) \\
Y'_2 &=& \overline{Y_2 \setminus (Y_2 \cap \nu_3)} &,& Y'_1 &=&
\overline{Y_1 \setminus (Y_1 \cap (\nu_2 \cup \nu_3))} \\
\tY'_2 &=& \rho^{-1} (Y'_2) &,& \tY'_1 &=& \rho^{-1} (Y'_1) \\
Z'(e) &=& \overline{Z(e) \setminus (Z(e) \cap \nu_3)} &,& \pa Z'(e) &=& Z(e)
\cap \pa \nu_3
\end{array}$$
Note that $\pa Y'_2 = \pa \nu_3$. Also, $\tY_2'\subset \tY_2$ is a
homotopy equivalence, as are all the other similarly related primed
and unprimed pairs. We denote by $\pa_{\rm ver}\nu'_2$ the closure
of $\pa\nu_2'\setminus \pa_{\rm hor}\nu_2'$. Similarly, we define
$\pa_{\rm ver}\widetilde \nu_2'$.

The pair $(\tY'_2,\tY'_1)$ is excisively equivalent to the pair $
(\widetilde{\nu}'_2, \pa \widetilde{\nu}'_2)$.  Analogously
$(Y'_2,Y'_1)$ is excisively equivalent to  $(\nu'_2, \pa \nu'_2)$.
We have the long exact sequence of the pairs
$$\begin{array}{ccccccccc}\label{lesHstarY2}
H_{*+1}(\widetilde{\nu}'_2, \pa_{\rm ver} \widetilde{\nu}'_2) & \to
& H_*(\tY'_1) & \to & H_*(\tY'_2) & \to &
H_*(\widetilde{\nu}'_2, \pa_{\rm ver} \widetilde{\nu}'_2) & \to & H_{*-1}(\tY'_1) \\
\downarrow & &  || & & \downarrow & & \downarrow & & || \\
H_{*+1}({\nu}'_2, \pa_{\rm ver} {\nu}'_2) & \to & H_*(Y'_1) & \to &
H_*(Y'_2) &
\to & H_*({\nu}'_2, \pa_{\rm ver} {\nu}'_2) & \to & H_{*-1}(Y'_1) \\
\downarrow & & & & & &  & &   \\
 0 & & & & & &  & &
\end{array}$$

Our descriptions of $\widetilde{\nu}_2$ and $\nu_2$ lead
to the following:
\begin{claim}
\begin{eqnarray*}
H_*(\widetilde{\nu}'_2,\pa_{\rm ver} \widetilde{\nu}'_2) & = &
\bigoplus_e ( H_{*-2}(Z'(e)) \otimes \scrA^*_{\ell^*(e)} ) \oplus
\bigoplus_e ( H_{*-4}(Z'(e)) \otimes
\Zee ) \\
H_*({\nu}'_2,\pa_{\rm ver} {\nu}'_2) & = & \bigoplus_e (
H_{*-2}(Z'(e)) \otimes Q_e ) \oplus \bigoplus_e ( H_{*-4}(Z'(e))
\otimes \Zee )
\end{eqnarray*}
and the projection map
$$\rho_* \colon H_*(\widetilde{\nu}'_2, \pa_{\rm ver} \widetilde{\nu}'_2) \to
H_*(\nu'_2, \pa_{\rm ver} \nu'_2)$$ is the natural one.
\end{claim}
Thus, $\rho_*$ is surjective with kernel $\bigoplus_e H_{*-2}(Z'(e)) \otimes
\scrA_{\ell^*(e)}$.

An elementary diagram chase in the long exact sequence above shows
that $H_2(\tY'_2) \to H_2(Y'_2)$ is surjective and its kernel is
identified with the kernel of
$$H_2(\widetilde{\nu}'_2, \pa_{\rm ver} \widetilde{\nu}'_2) \to H_2({\nu}'_2,
\pa_{\rm ver} {\nu}'_2) . $$

Thus,  we have established the following relationship between $H_*(\tY_2)$ and
$H_*(Y_2)$:
\begin{lemma}\label{tY2Y2} There is an exact sequence
$$0 \to \bigoplus_e H_{*-2}(Z(e)) \otimes \scrA_{\ell^*(e)} \to
H_*(\tY_2) \to H_*(Y_2) \to 0 .$$ The kernel is represented by the
fundamental classes of the $\Pee(\ell)$-bundles over cycle
representatives for the classes in $H_{*-2}(Z(e))$.
\end{lemma}

Exactly the same argument can be applied to
$\pa\widetilde{\nu}_3\to \pa {\nu}_3$. Recall from the definition
that $\pa Z'(e)$ is contained in $\pa\nu_3$. For any edge $e$ and
any two-face $f$ containing $e$ in its closure denote by
$\pa_fZ'(e)$ the components of $\pa Z(e)$ given by $\pa\nu(f)\cap
Z'(e)$, so that $\pa Z'(e)=\coprod_{f| e\prec f}\pa_fZ'(e)$.

\begin{lemma}\label{local2} For any two-face $f$ of $\partial \Delta$
there is a short exact sequence
$$0 \to \bigoplus_{e\in E;\ e \prec f}
H_{*-2}(\pa_fZ'(e)) \otimes \scrA_{\ell^*(e)} \to H_*(\pa
\widetilde{\nu}(f)) \to H_*(\pa \nu(f)) \to 0 .$$ The kernel is
represented by the fundamental classes of the
$\Pee^1(\ell)$-bundles over the cycles in $\pa_fZ'(e)$ for the
classes in $H_{*-2}(\pa_fZ'(e))$.
\end{lemma}

In order to carry out the next step in the comparison, we need to
understand the relative version, that is to say the map
$H_*(\tY'_2, \pa\widetilde{\nu}_3) \to H_*(Y'_2, \pa\nu_3)$.

\begin{prop}\label{3.4}
$H_*(\tY'_2, \pa \widetilde{\nu}_3) \to H_*(Y'_2, \pa \nu_3)$ is an
isomorphism for $* \leq 2$ and is surjective with kernel
$\bigoplus_{e\in E} H_1(Z'(e), \pa Z'(e)) \otimes \scrA_{\ell^*(e)}$
for $* = 3$.
\end{prop}

\begin{proof}
We consider $\pa \widetilde{\nu}_3 \subset (\pa \widetilde{\nu}_3
\cup \widetilde{\nu}'_2) \subset \tY'_2$ and $\pa {\nu}_3 \subset
(\pa \nu_3 \cup {\nu}'_2) \subset Y'_2$.  Of course, $(\tY'_2, \pa
\widetilde{\nu}_3 \cup \widetilde{\nu}'_2) \simeq (\tY'_1, \pa
\tY'_1)$ and $(Y'_2, \pa {\nu}_3 \cup {\nu}'_2) \simeq (Y'_1, \pa
Y'_1)$.  Thus we have
$$\hspace{-.5in} \begin{array}{ccccccccc} H_{*+1}(\tY'_1, \pa
\tY'_1) & \to & H_*(\widetilde{\nu}'_2, \pa_{\rm hor}
\widetilde{\nu}'_2) & \to & H_*(\tY'_2, \pa \widetilde{\nu}_3) & \to
& H_*(\tY'_1, \pa \tY'_1) & \to & H_{*-1}(\widetilde{\nu}'_2,
\pa_{\rm hor} \widetilde{\nu}'_2) \\
|| & & \downarrow & & \downarrow & & ||
& & \downarrow  \\
H_{*+1}(Y'_1, \pa Y'_1) & \to & H_*({\nu}'_2, \pa_{\rm hor}
{\nu}'_2) & \to & H_*(Y'_2, \pa {\nu}_3) & \to & H_*(Y'_1, \pa Y'_1)
& \to & H_{*-1}({\nu}'_2, \pa_{\rm hor} {\nu}'_2) .
\end{array}$$
Again, the description of $\widetilde{\nu}_2$ and $\nu_2$ leads
immediately to the following:
\begin{claim}
\begin{eqnarray*}
H_*(\widetilde{\nu}'_2, \pa_{\rm hor} \widetilde{\nu}'_2) & = &
\bigoplus_{e\in E} ( H_{*-2}(Z'(e), \pa Z'(e)) \otimes
\scrA_{\ell^*(e)} ) \oplus \bigoplus_{e\in E}
H_*(Z'(e), \pa Z'(e)) \\
H_*(\nu'_2, \pa_{\rm hor} \nu'_2) & = & \bigoplus_{e\in E}
H_*(Z'(e), \pa Z'(e))
\end{eqnarray*}
and the projection map induces the natural one on homology.
\end{claim}
Thus, $H_*(\widetilde{\nu}'_2, \pa_{\rm hor} \widetilde{\nu}'_2)
\to H_*(\nu'_2, \pa_{\rm hor} \nu'_2)$ is surjective with kernel
$H_{*-2}(Z'(e), \pa Z'(e))$. In particular,
$H_*(\widetilde{\nu}'_2, \pa_{\rm hor} \widetilde{\nu}'_2) \to
H_*(\nu'_2, \pa_{\rm hor} \nu'_2)$ is an isomorphism for $* \leq
2$. Plugging this into the commutative diagram, the five lemma
tells us that $H_*(\tY'_2, \pa \widetilde{\nu}_3) \to H_*(Y'_2,
\pa \nu_3)$ is an isomorphism for $* \leq 2$.

Now we consider $* = 3$.
$$\hspace{-.35in} \begin{array}{ccccccccc} H_{4}(\tY'_1, \pa
\tY'_1) & \to & H_3(\widetilde{\nu}'_2, \pa_{\rm hor}
\widetilde{\nu}'_2) & \to & H_3(\tY'_2, \pa \widetilde{\nu}_3) &
\to & H_3(\tY'_1, \pa \tY'_1) & \to & H_{2}(\widetilde{\nu}'_2,
\pa_{\rm hor} \widetilde{\nu}'_2) \\
\, \downarrow = & & \downarrow & & \downarrow & & \, \downarrow =
& & \, \downarrow \simeq  \\
H_{4}(Y'_1, \pa Y'_1) & \to & 0 & \to & H_3(Y'_2, \pa {\nu}_3) & \to
& H_3(Y'_1, \pa Y'_1) & \to & H_{2}({\nu}'_2, \pa_{\rm hor}
{\nu}'_2) .
\end{array}$$
A standard diagram chase shows that $H_3(\tY'_2, \pa
\widetilde{\nu}_3) \to H_3(Y'_2, \pa \nu_3)$ is onto, with kernel
the image of $H_3(\widetilde{\nu}'_2, \pa_{\rm hor}
\widetilde{\nu}'_2) \to H_3(\tY'_2, \pa \widetilde{\nu}_3)$.  Of
course,
$$H_3(\widetilde{\nu}'_2, \pa_{\rm hor} \widetilde{\nu}'_2) = \bigoplus_{e\in E}
H_1(Z'(e), \pa Z'(e)) \otimes \scrA_{\ell^*(e)} .$$  To complete the
argument, we need to know that $H_4(\tY'_1, \pa \tY'_1) \to
H_3(\widetilde{\nu}'_2, \pa_{\rm hor} \widetilde{\nu}'_2)$ is the
trivial map. That is the content of the next lemma.

\begin{lemma} \label{Homologynutilde2}
\begin{equation}\label{H4} H_4(\tY'_1, \pa\tY'_1) \to
H_3(\widetilde{\nu}'_2, \pa_{\rm hor} \widetilde{\nu}'_2)
\end{equation} is the trivial map.
\end{lemma}
\begin{proof}
The morphism factors as
$$\hspace{-.15in} H_4(\tY'_1, \pa \tY'_1) \to H_3(\pa \tY'_1) \to H_3(\pa \tY'_1,
\widetilde{\nu}_3 \cap \pa \tY'_1) \simeq H_3(\pa_{\rm ver}
\widetilde{\nu}'_2, \pa_{\rm hor} \pa_{\rm ver} \widetilde{\nu}'_2)
\to H_3(\widetilde{\nu}'_2, \pa_{\rm hor} \widetilde{\nu}'_2) \ ,$$
where we introduce the notation \begin{eqnarray*} \pa_{\rm ver}
\nu'_2 & = & \nu'_2 \cap Y'_1 \\ \pa_{\rm hor}\pa_{\rm ver} \nu'_2 &
= & \pa_{\rm ver}\pa_{\rm hor} \nu'_2 = \pa (\nu'_2 \cap \pa \nu_3)
= \pa_{\rm ver} \nu'_2 \cap \pa \nu_3 \ .
\end{eqnarray*}
The isomorphism is from excision.  The pair $(\widetilde{\nu}'_2,
\pa_{\rm hor} \widetilde{\nu}'_2)$ is a locally trivial fiber
bundle over $\coprod_e(Z'(e),\pa Z'(e))$ with fiber over any $x\in
Z'(e)$ a regular neighborhood of $\scrA_{\ell^*(e)}$. Hence,
$H_3(\widetilde{\nu}'_2, \pa_{\rm hor} \widetilde{\nu}'_2))\simeq
\bigoplus_eH_1(Z'(e),\pa Z'(e))\otimes \scrA_{\ell^*(e)}$ is free
abelian. On the other hand, $(\pa_{\rm ver} \widetilde{\nu}'_2,
\pa_{\rm hor} \pa_{\rm ver} \widetilde{\nu}'_2)$ is a locally
trivial fiber bundle over $\coprod_e(Z'(e),\pa Z'(e))$ with fiber
over any $x \in Z'(e)$ identified with the boundary of a regular
neighborhood of the cone point $0_e$ in $V(e)$. The latter is the
lens space $L(n(e), -1)$. Thus, $H_3(\pa_{\rm ver}
\widetilde{\nu}'_2, \pa_{\rm hor} \pa_{\rm ver}
\widetilde{\nu}'_2) = \bigoplus_e H_2(Z'(e), \pa Z'(e)) \otimes
\Tor_{\ell^*(e)}$ is a torsion group. Hence, the morphism in
Equation~(\ref{H4}) is the zero map.
\end{proof}

This completes the proof of Proposition~\ref{3.4}.
\end{proof}

Note that the three inclusions of
\begin{eqnarray*}
(\tY'_2, \pa \widetilde\nu_3) & \hookrightarrow & (\tY,
\widetilde\nu_3) \\
(Y'_2, \pa \nu_3) & \hookrightarrow & (Y_2, \nu_3) \\
(Z'(e), \pa Z'(e)) & \hookrightarrow & (Z(e), Z(e) \cap \nu_3)
\end{eqnarray*}
each consist of excisive pairs.

\subsection{The comparison of $\tY$ and $Y$}
In the previous subsection, to study the relationship between
$\tY_2$ and $Y_2$, for each edge $e$ of $\pa \Delta$ we described
$\scrS(e)$ in terms of its combinatorial model which is the cell
complex dual to $\scrT|_e$. Analogously, in order to understand
the relationship between $\tY$ and $Y$ we need to understand the
homology of the surfaces $\scrS(f)$ for two-faces $f$ of $\pa
\Delta$. Of course by Proposition~\ref{propmodels} we have
analogous result for each two-face $f$: the combinatorial model
for $\scrS(f)$ is the cell complex dual to the triangulation of
$\stackrel{\circ}{f}$. We study $\scrS(f)$ using this cell
complex.

Let $\Sigma$ be a finite, connected cell complex dual to the triangulation of
the interior of a two-face. Then every cell of $\Sigma$ has dimension $\le 2$,
 every vertex $v$ has valence, denoted $o(v)$, which is $\le 3$, and
every edge has two vertices. Let $\scrS$ be a compact algebraic
variety of dimension $\le 2$ such that $\Sigma$ is a combinatorial
model for $\scrS$. It follows from the definition that for every
edge $\alpha$ of $\Sigma$ the algebraic subset $\overline
S_\alpha$ of $\scrS $ is isomorphic to $\Pee^1$. We also denote
this subset by $\Pee^1(\alpha)$.

Define $\Sigma_2\subset \Sigma$ to be the sub-cell complex which
is the closure of the union of the two-cells of $\Sigma$. Let
$\Sigma_1\subset \Sigma$ be the sub-cell complex consisting the
closure of the union of all edges of $\Sigma$ that are not
contained in $\Sigma_2$.

\begin{figure}[ht] 
\begin{center}
$\begin{array}{c@{\hspace{.5in}}c}
\includegraphics[width=2in]{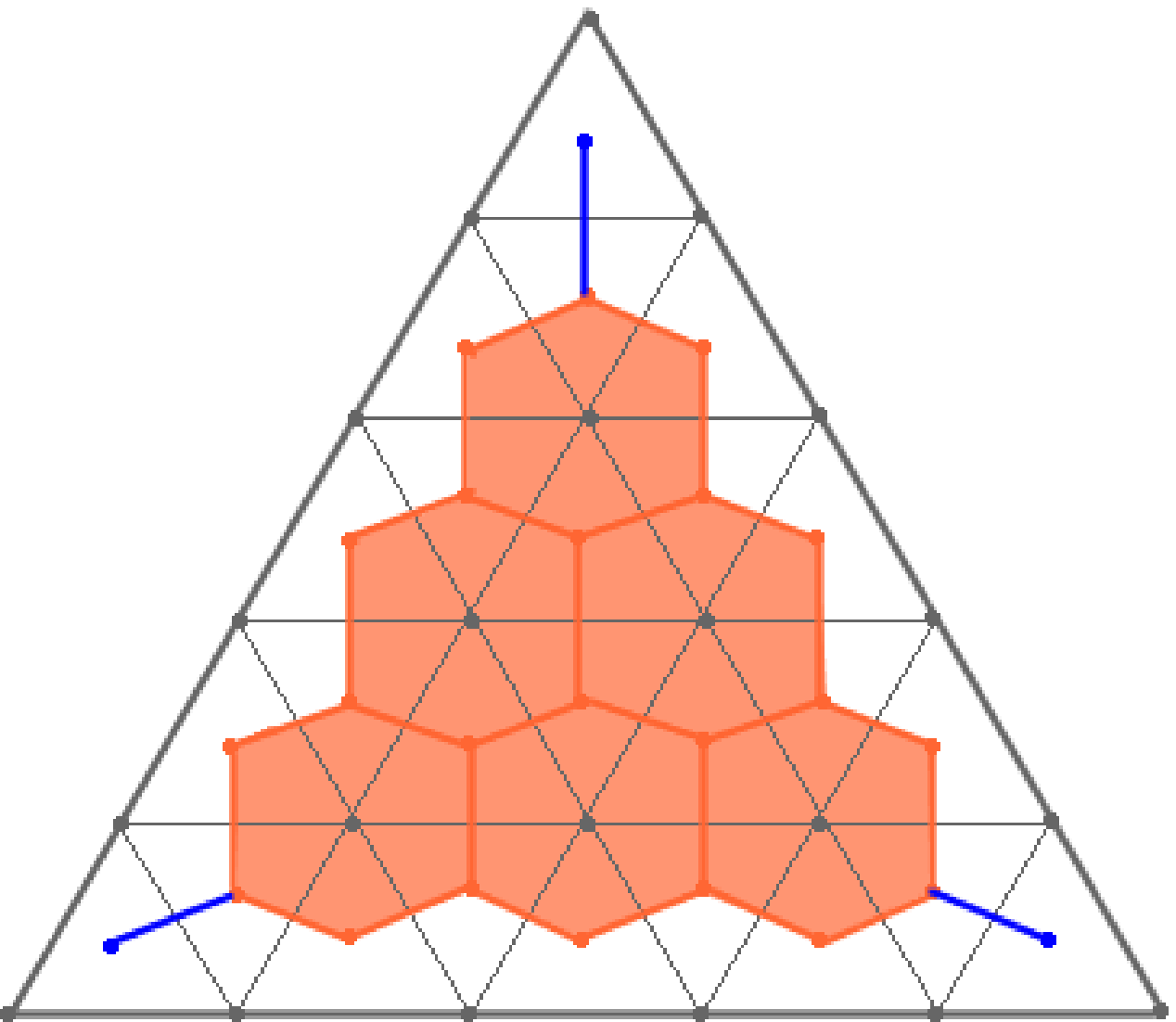} &
\includegraphics[width=2in]{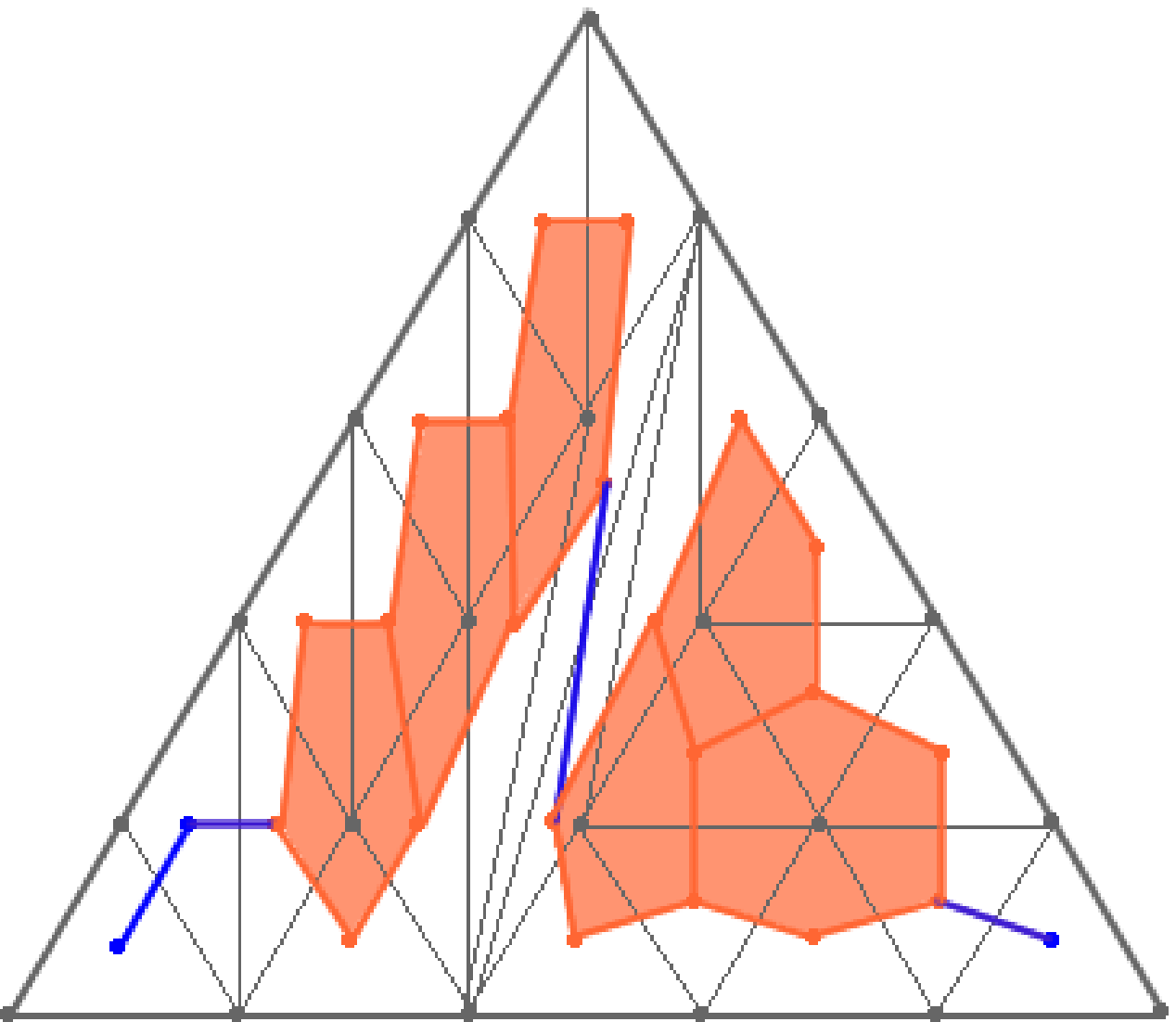}
\end{array}$
\end{center}
\caption{$\Sigma_1$ and $\Sigma_2$ for the triangulations of
Figure \protect\ref{figuretriang}}
\label{figuretriangshaded}
\end{figure}

Notice that $\Sigma_2$ is a disjoint union of
contractible homogeneously two-dimensional complexes, and
$\Sigma_1$ is a disjoint union of trees. Furthermore,
$\Sigma_2\cap \Sigma_1$ consists of a set of trivalent vertices of
$\Sigma$. If $A$ is a component of $\Sigma_2$ we denote by $\scrS
_A$ the corresponding, possibly reducible, surface in $\scrS$. If
$X$ is a component of $\Sigma_1$ we denote by $\scrS _X$ the
corresponding curve (which is a union of irreducible rational
components) in $\scrS $. We denote by $\scrC_2$ the set of
two-cells in $\Sigma$, by $\scrE_2$ the set of edges in
$\Sigma_2$, and by $\scrB_2$  the set of interior edges of
$\Sigma_2$. Lastly, we denote by $\scrE_1$ the set of edges in
$\Sigma_1$.

The following is immediate from the Meyer-Vietoris sequence.

\begin{lemma}\label{homology}
Let $\Sigma$ be a combinatorial model for $\scrS $ as above. Then
$$\begin{array}{rcl}
H_0(\scrS) &=& \Zee \\
H_1(\scrS) &=& \bigoplus_{A \in \pi_0(\Sigma_2)}H_1(\scrS_A) \\
 H_2(\scrS) &=& \bigoplus_{A \in
\pi_0(\Sigma_2)}H_2(\scrS_A)\oplus \bigoplus_{\alpha \in \scrE_1}H_2(\Pee^1(\alpha)) \\
H_3(\scrS) &=& \bigoplus_{A \in \pi_0(\Sigma_2)}H_3(\scrS_A) \\
H_4(\scrS) &=& \bigoplus_{A \in \pi_0(\Sigma_2)}
H_4(\scrS_A).\end{array}$$
\end{lemma}

Recall that
$$\scrS(f)= \rho^{-1}(0_f).$$
We now wish to compute the cohomology of
$(\Pee_{\scrF'_f(N_f)},\Pee_{\scrF'_f(N_f)}\setminus \scrS(f))$.
 Of course, by Proposition \ref{propmodels},  $\Sigma(\scrF'_f)$ is a
 combinatorial model
for $\scrS(f)$. Clearly, $\Sigma(\scrF'_f)$ is identified with a
polyhedral subset of the face $f$ of $\Delta$ and hence is
embedded in the plane spanned by $f$.
 To
simplify notation we set $\Sigma=\Sigma(\scrF'_f)$ and
$S=\scrS(f)$.

 Consider the pair
$$\left(\Pee_{\scrF'_f(N_f)},\Pee_{\scrF'_f(N_f)}\setminus S\right).$$
Since $\Pee_{\scrF'_f(N_f)}$ deformation retracts to $S$, by
Lefschetz duality, we have
$$H^i\left(\Pee_{\scrF'_f(N_f)},\Pee_{\scrF'_f(N_f)}\setminus S\right)=H_{6-i}\left( S \right).$$

The following lemma computes these groups:

\begin{lemma} \label{lemma:cohomgrps}
For any $f\in F$, we have
\begin{eqnarray*}
H_6( S ) & = & 0\\
H_4( S ) & = & \bigoplus_{\left\{c\in \scrC_2\right\}} \Zee[\overline S_c]  \\
H_0( S ) & = & \Zee,
\end{eqnarray*}
and all odd homology groups vanish. Furthermore, $H_2( S )$ is
free abelian with rank
$$b_2( S )= \ell^*(f)+\sum_{\left\{{\rm edges}\ e\left|\right.e\prec
f\right\}}\ell^*(e) + v_f - 3. \label{eqn:b2S}$$
\end{lemma}

\begin{proof}
Let $\scrT_f$ be the induced triangulation of $f$. We say that an
edge of a two-dimensional triangulation of a surface is {\em free}
if it is the face of exactly one two-simplex.  There are three types
of simplifications we wish to perform on triangulations of surfaces:
\begin{itemize}
\item[Type A] Remove a two-simplex whose closure meets the
boundary in exactly one edge
 and remove that free edge as well.
\item[Type B] Remove a two-simplex with exactly two free edges and
remove both free edges and the vertex they have in common.
\item[Type C] Remove a two-simplex with three free edges, all its
edges, and all its vertices.
\end{itemize}
A triangulated surface with boundary is said to be {\em shellable}
if there is a finite sequence of these simplifications that
reduces it to the empty triangulation. Notice that each of these
operations replaces a triangulated surface by a triangulated
subsurface with one fewer two-simplex. As we perform a sequence of
these operations we obtain a decreasing sequence of triangulations
on subsurfaces of $f$. Dual to this sequence of triangulations is
the sequence of dual cell complexes. Each time we remove simplices
from one of the triangulations, we remove their dual cells from
the dual cell complex. Thus, we produce a decreasing sequence of
subcell complexes $\Sigma=\Gamma_0\supset \Gamma_1\supset
\cdots\supset \Gamma_m=\emptyset$. Corresponding to the sequence
of subcell complexes there is a decreasing sequence of
$T_{N_f}$-invariant locally closed algebraic subsets
$S=\scrS_0\supset \scrS_1\supset \cdots \supset
\scrS_m=\emptyset$: As we remove a cell of $\Gamma_i$ we remove
the corresponding $T_{N_f}$ orbit from $\scrS_i$. Suppose that
$\scrS_{i+1}$ is obtained from $\scrS_i$ by an operation of Type
A. Then $\scrS_{i+1}$ is obtained by removing  a contractible
subset of a closed irreducible (complex) surface component of
$\scrS_i$. The subset is disjoint from all other irreducible
components of $\scrS_i$.  Hence,
$b_4(\scrS_{i+1})=b_4(\scrS_i)-1$, and all other homology groups
are the same. Suppose instead that $\scrS_{i+1}$ is obtained from
$\scrS_i$ by an operation of Type B. Then $\scrS_{i+1}$ is
obtained by removing  a contractible subset of an irreducible
component  of $\scrS_i$ isomorphic to $\Pee^1$. Again, the subset
is disjoint from all other irreducible components of $\scrS_i$.
Hence, $b_2(\scrS_{i+1})=b_2(\scrS_i)-1$, and all other homology
groups are the same. Lastly, suppose that $\scrS_{i+1}$ is
obtained from $\scrS_i$ by an operation of Type C. Then
$\scrS_{i+1}$ is obtained by removing  an isolated point of
$\scrS_i$. Hence, $b_0(\scrS_{i+1})=b_0(\scrS_i)-1$, and all other
homology groups are the same.

\begin{claim}
Let $\tau$ be a compact triangulated contractible surface. Then
there is a sequence of operations of Types A,B, and C that reduce
$\tau$ to the empty surface.
\end{claim}

\begin{proof}
The argument is by induction on the number of two-simplices in $\tau$. If there
is a two-simplex of $\tau$ whose intersection with the boundary is an edge,
then we can perform an operation of Type A on this two-simplex producing a
contractible surface with one fewer two-simplex. Suppose there are no such
two-simplices. Then every two-simplex that contains an edge of the boundary has
all its vertices in the boundary. Orient the boundary and let $\sigma$ be a
two-simplex containing an edge $\alpha$ of the boundary with third vertex
$p_\alpha$. Define the distance from $\alpha$ to $p_\alpha$ along the boundary
to be the number of edges moving in the positive direction that separate
$\alpha$ from $p_\alpha$. If this distance is $1$, then $\sigma$ has at least
two edges in the boundary. If it has exactly two, then we can perform an
operation of Type B on $\sigma$. If it has all three edges in the boundary,
then we can perform an operation of Type C on $\sigma$. If the distance along
the boundary from $\alpha$ to $p_\alpha$ is greater than one, let $\alpha'$ be
the edge next to $\alpha$ in positive direction along the boundary. Clearly,
the distance along the boundary from $\alpha'$ to $p_{\alpha'}$ is less than
that from $\alpha$ to $p_\alpha$. Eventually we arrive at an  edge $\alpha_0$
with the distance from $\alpha_0$ to $p_{\alpha_0}$ being $1$, completing the
induction.
\end{proof}

Given this,  it follows from the above description and induction
that all the homology groups of $S$ are free abelian and all those
except $H_4(S)$, $H_2(S)$, and $H_0(S)$ vanish. The groups $H_4$ and
$H_0$ are clearly as stated. The rank of $H_2$ is the number of
operations of Type  B. A direct combinatorial argument shows that
the number of operations of Type B is as given in the statement.
\end{proof}

Since $(\tY,\tY'_2)$ is excisively equivalent to
$(\widetilde{\nu}_3,\pa \widetilde{\nu}_3)$, we have:
\begin{corollary}\label{nu3}
$$H_i(\widetilde{\nu}_3,\pa\widetilde{\nu}_3) \simeq H_i(\tY,\tY'_2)
= \bigoplus_{f\in F}H_0(Y\cap \scrO_f)\otimes
H_i\left(\Pee_{\scrF'_f(N_f)},\Pee_{\scrF'_f(N_f)}\setminus
\scrS(f)\right).$$ Of course, since $\widetilde{\nu}_3$ is
homotopy equivalent to $\scrS(f)$, we see that
$H_{odd}(\widetilde{\nu}_3) = 0$ and $H_{ev}(\widetilde{\nu}_3)$
is free abelian.
\end{corollary}

\subsection{The computation of $\rho_*\colon H_1(\tY)\to H_1(Y)$}

\begin{prop}\label{proprho2skel}
The morphism
$$\rho_* \colon H_1(\tY) \to H_1(Y)$$
is an isomorphism, and $H_1(\tY)$ is identified with
the quotient of $N$ by the lattice spanned by the
2-skeleton of $\pa \Delta$.
\end{prop}

\begin{proof}
By Lemma~\ref{tY2Y2}  the restriction of $\rho$ induces an isomorphism
$$H_1(\tY_2) \stackrel{\cong}{\longrightarrow} H_1(Y_2) .$$
We have the exact sequence
$$\begin{array}{ccccccccccc}
\cdots & \to & H_2(\tY) & \to & H_2(\tY,\tY_2) & \to & H_1(\tY_2)
& \to &
H_1(\tY) & \to & 0 \\
 &  & \downarrow & & \downarrow & & \, \downarrow\cong  & & \downarrow & & || \\
\cdots & \to & H_2(Y) & \to & H_2(Y,Y_2) & \to & H_1(Y_2) & \to &
H_1(Y) & \to  & \ \ 0 \ . \end{array} $$ To prove that the map
$\rho_*\colon H_1(\tY)\to H_1(Y)$ is an isomorphism we need only
see that the natural map $H_2(\tY,\tY_2)\to H_2(Y,Y_2)$ is onto.
By excision, it is equivalent to show that $H_2(\widetilde{\nu}_3,
\pa \widetilde{\nu}_3) \to H_2(\nu_3, \pa\nu_3)$ is onto.

\begin{claim} \label{claimnu3onto}
$H_2(\widetilde{\nu}_3, \pa \widetilde{\nu}_3) \to H_2(\nu_3, \pa
\nu_3)$ is onto and the range is torsion.
\end{claim}
\begin{proof}
Since $\nu_3$ is a disjoint union of cones, $H_2(\nu_3, \pa\nu_3)
\cong H_1(\pa \nu_3)$.

Lemma~\ref{local2} shows that $H_1(\pa \widetilde{\nu}_3)\to
H_1(\pa \nu_3)$ is onto. Thus, to complete the proof of
surjectivity, we need only see that $H_2(\widetilde{\nu}_3, \pa
\widetilde{\nu}_3)\to H_1(\pa \widetilde{\nu}_3)$ is surjective,
but that follows because $H_1(\widetilde{\nu}_3)=0$ by
Corollary~\ref{nu3}.

To establish that $H_2(\nu_3, \pa\nu_3)$ is torsion, we must show
that for each two-face $f$, $H_1(\pa \nu_3(f)) = 0$.  Of course,
$\pa\nu_3(f)$ is homotopy-equivalent to $V(f)^* = V(f) \setminus
0_f$.  But $V(f)$ is a three-dimensional toric variety and
$V(f)^*$ is the union of all the torus orbits of codimension $<
3$, i.e., $V(f)^* = V(f)_2$ (with notation analogous to that for
the $X_i \subset X$).  The result now follows from the
observations that $V(f)_0 = T_{N_f}$, that $$H_1(V(f)_1) =
\frac{N_f}{\Span V(f)}$$ (which is torsion) by the exact sequence
for the pair $V(f)_0 \hookrightarrow V(f)_1$, and that
$$H_1(V(f)_2) = \frac{N_f}{\langle L \cap f^{(1)} \rangle}$$
by the exact sequence for the pair $V(f)_1 \hookrightarrow
V(f)_2$.
\end{proof}

This completes the proof of Proposition~\ref{proprho2skel}
\end{proof}

\subsection{The computation of $\rho_* \colon H_2(\tY) \to H_2(Y)$} \label{subsecCompH2}

In Corollary~\ref{2.23} we showed that in the hypersurface case
$$\rk H_2(\tY) = \# V + \sum_e \ell^*(e) + \sum_f
\ell^*(f) (1 + \ell^*(f^\vee)) - 4 $$ and in the complete
intersection case that
$$\rk H_2(\tY) = \# V + \sum_e \ell^*(e) + \sum_f
\ell^*(f) (1 +\sum_{J \subset I;\ J\not=\emptyset} (-1)^{3 - |J|}
\ell^*( \sum_{j \in J} f_j^\vee)-n.$$ We also gave  a formula in
terms of the lattice and $\Delta$ for the torsion part of this
group. The new result here is:

\begin{prop} \label{propH2tY}  The morphism
$$\rho_* \colon H_2(\tY) \to H_2(Y)$$ is onto.
%and $$\rk H_2(\tY) = \sum_f \ell^*(f) \cdot \# (\scrO_f \cap Y) +
%\sum_e \ell^*(e) + \# V - n$$
\end{prop}
\begin{proof}
We begin by computing with the two excisive triples
\begin{eqnarray*}
\tY'_2 \subset \tY \to (\tY, \tY'_2) & {\simeq} &
(\widetilde{\nu}_3, \pa \widetilde{\nu}_3) \\
Y_2' \subset Y \to (Y, Y'_2) & {\simeq} & ({\nu}_3, \pa {\nu}_3)
\end{eqnarray*}
from which we have:
$$\begin{array}{ccccccccc}
 & & 0          & & 0          & & 0          & &  \\
 & & \downarrow & & \downarrow & & \downarrow & &   \\
0 & & \oplus_e \scrA_{\ell^*(e)} & \to & \Ker_2 & \to & \Ker'_2
& &  \\
|| & & \downarrow & & \downarrow & & \downarrow & &  \\
H_{3}(\widetilde{\nu}_3, \pa \widetilde{\nu}_3) & \to &
H_2(\tY'_2) & \to & H_2(\tY) & \to &
H_2(\widetilde{\nu}_3, \pa \widetilde{\nu}_3) & \to & H_{1}(\tY'_2) \\
\downarrow & & \downarrow & & \downarrow & &  \downarrow
& & \, \downarrow =  \\
H_{3}({\nu}_3, \pa {\nu}_3) & \to & H_2(Y'_2) & \to & H_2(Y) & \to
& H_2({\nu}_3, \pa {\nu}_3) & \to & H_{1}(Y'_2) \\
 & & \downarrow & &   & & \downarrow & & \\
 & & 0          & &   & & 0          & &
\end{array}$$
where the exactness of the second and fourth columns was
established in Lemmas~\ref{tY2Y2} and ~\ref{local2} respectively,
and we take by definition $\Ker_2$ to be the kernel of $H_2(\tY)
\to H_2(Y)$. Here Corollary~\ref{nu3} implies that
$H_3(\widetilde{\nu}_3, \pa \widetilde{\nu}_3) = 0$. By
Lemma~\ref{tY2Y2},     the sequence
$$0 \to \scrA_{\ell^*(e)} \to H_2(\tY'_2) \to H_2(Y'_2) \to 0$$
is exact.  {\bf N.B.:} We are not claiming that the top row is
exact.

By Claim~\ref{claimnu3onto}, $H_2(\widetilde{\nu}_3, \pa
\widetilde{\nu}_3)$ is onto. Finally, $H_1(\tY'_2) = H_1(Y'_2)$ by
Lemma~\ref{tY2Y2}.
 A straightforward
diagram chase completes the proof of Proposition~\ref{propH2tY}.
\end{proof}

\subsection{The computation of $\rho_* \colon H_3(\tY) \to H_3(Y)$}

For each edge $e$ of $\partial \Delta$ we have the (incomplete)
algebraic curve $Z(e)=Y\cap \scrO_e$. Its completion $\widehat
Z(e)= Y\cap \overline{\scrO_e}$ is a compact complex curve. Since
$\overline{\scrO_e}$ has singularities at most at isolated points,
$\widehat Z(e)$ is a smooth riemann surface. Recall that for each
edge $e$ of $\pa\Delta$ and each lattice point $\ell\in
\stackrel{\circ}{e}$ we have a ruled surface $R(\ell)$ over
$\widehat Z(e)$. The generic fiber of $R(\ell)\to \widehat Z(e)$
is a smooth rational curve denoted $\Pee^1(\ell)$. We also have
$H_3(R(\ell))=H_1(\widehat Z(e))\otimes H_2(\Pee^1(\ell))$. Taking
the union over all $\ell\in\stackrel{\circ}{e}$ gives a surface
whose $H_3$ is identified with $H_1(\widehat Z(e))\otimes
\scrA_{\ell^*(e)}$.

 Here we shall prove:
\begin{prop} \begin{enumerate} \item The kernel of $\rho_* \colon H_3(\tY) \to
H_3(Y)$ is identified with $\bigoplus_e H_1(\widehat{Z}(e)) \otimes
\scrA_{\ell^*(e)}$ given by the inclusion of the ruled surfaces
$\coprod_{e; \ell\in\stackrel{\circ}{e}} R(\ell)$
$$\bigoplus_{e;\ell\in \stackrel{\circ}{e}} H_3(R(\ell)) \to H_3(\tY) .$$
\item The image of $\rho_* \colon H_3(\tY) \to H_3(Y)$ is
$\Image(H_3(Y_2) \to H_3(Y))$.
\end{enumerate}
\end{prop}

\begin{proof}
First let us consider the kernel of the map $\rho_* \colon H_3(\tY)\to H_3(Y)$.
Consider the long exact commutative diagram associated to the pairs
\begin{eqnarray*}
\widetilde{\nu}_3 \subset \tY \to (\tY, \widetilde{\nu}_3) &
{\simeq} &
(\tY'_2, \pa \widetilde{\nu}_3) \\
\nu_3 \subset Y \to (Y, \nu_3) & {\simeq} & (Y'_2, \pa {\nu}_3).
\end{eqnarray*}
We obtain
$$\hspace{-.25in} \begin{array}{ccccccccc}
 & &   &     & 0 & & 0 & & \\
 & &   &     & \downarrow & & \downarrow & & \\
 & & 0 & \to & \Ker_3 & \to & \oplus_e H_1(Z'(e), \pa Z'(e)) \otimes \scrA_{\ell^*(e)} & &  \\
 & & || & & & & \downarrow & &  \\
H_{4}(\tY'_2, \pa \widetilde{\nu}_3) & \to &
H_3(\widetilde{\nu}_3) & \to & H_3(\tY) & \to & H_3(\tY'_2, \pa
\widetilde{\nu}_3) & \to &
H_2(\widetilde{\nu}_3) \\
 & & \downarrow & & \downarrow & &  \downarrow
& & \downarrow   \\
& & H_3(\nu_3) & \to & H_3(Y) & \to
& H_3(Y'_2, \pa {\nu}_3) & \to & H_{2}(\nu_3) \\
& & ||         &     &        &
& \downarrow            &     & || \\
& & 0          &     &        & & 0                     &     & \
\  0 \ .
\end{array}$$
The short exact sequence in the fourth column is established in
Proposition~\ref{3.4}. Corollary~\ref{nu3} shows that
$H_3(\widetilde \nu_3)=0$. Of course, since $\nu_3$ is a disjoint
union of cones, its third homology vanishes.

 Now a diagram
chase shows that $\Ker_3 = \Ker (\rho_* \colon H_3(\tY) \to H_3(Y))$ sits in
the sequence
\begin{equation}\label{Kseq} 0 \to \Ker_3 \to \oplus_e H_1(Z_e, \pa
Z_e) \otimes \scrA_{\ell^*(e)} \stackrel{\pa}{\to}
H_2(\widetilde{\nu}_3) .\end{equation} The following lemma
completes the computation of $\Ker_3$.

\begin{lemma}
The natural map
$$\coprod_{e;\ell\in \stackrel{\circ}{e}} R(\ell) \to \widetilde{Y}$$
induces an isomorphism from
$$\bigoplus_{e;\ell\in\stackrel{\circ}{e}} H_3(R(\ell))
\stackrel{\simeq}{\longrightarrow} \Ker_3 \ .$$
\end{lemma}

\begin{proof}
 The ``boundary morphism''
$\pa$ in Equation~\ref{Kseq} is the composition of natural
homomorphisms.  The first of these
\begin{equation}\label{firsthomo}\bigoplus_e H_1(Z'(e), \pa Z'(e)) \otimes \scrA_{\ell^*(e)}
\to \bigoplus_e H_0(\pa Z'(e)) \otimes
\scrA_{\ell^*(e)}\end{equation} is given by the direct sum of the
boundary homomorphism in relative homology on the first factor
tensored with the identity homomorphism on the $\scrA_{\ell^*(e)}$
term. The group
\begin{eqnarray*}
\bigoplus_e H_0(\pa Z'(e)) \otimes \scrA_{\ell^*(e)} & = &
\bigoplus_e
\bigoplus_{e \prec f} H_0(\scrO_f \cap Y) \otimes \scrA_{\ell^*(e)} \\
 & = & \bigoplus_f \bigoplus_{e \prec f} H_0(\scrO_f \cap Y) \otimes
\scrA_{\ell^*(e)} \\
  & = & \bigoplus_f  H_0(\scrO_f \cap Y) \otimes \bigoplus_{e \prec f}
  \scrA_{\ell^*(e)}
  \end{eqnarray*}
  where the first isomorphism follows from the fact that $\pa Z'(e)
  = \coprod_{e \prec f} (\scrO_f \cap Y) \times S^1$.  The
  second homomorphism
$$\bigoplus_f [ H_0(\scrO_f \cap Y) \otimes \bigoplus_{e \prec f}
  \scrA_{\ell^*(e)} ] \to H_2(\widetilde \nu_3)=\bigoplus_f [H_0(\scrO_f \cap Y) \otimes
  H_2(\scrS(f))]$$
is given by $\bigoplus_f \Id_{H_0(\scrO_f \cap Y)} \otimes \left[
\bigoplus_{e \prec f} \scrA_{\ell^*(e)} \to H_2(\scrS(f))\right]$.
The main point is to show the injectivity of this last morphism.
That is to say, for all $f$, the map $ \bigoplus_{e \prec f}
\scrA_{\ell^*(e)} \to H_2(\scrS(f))$ is injective. In light of the
isomorphism $H_2(\scrS(f)) \simeq H_2(\widetilde{\nu}_3(f))$, and
the fact that the morphism $\bigoplus_{e \prec f} \scrA_{\ell^*(e)}
\to H_2(\scrS(f))$ factors through the morphism $H_2(\pa
\widetilde{\nu}_3(f)) \to H_2(\scrS(f))$, which is itself injective
by the long exact sequence of the pair (the previous term being
$H_3(\widetilde{\nu}_3, \pa \widetilde{\nu}_3) = 0$), it suffices to
show that $\bigoplus_{e \prec f} \scrA_{\ell^*(e)} \to H_2(\pa
\widetilde{\nu}_3(f))$ is injective.  This last follows by
Lemma~\ref{local2}. Thus $K$ is the kernel of the  homomorphism in
Equation~\ref{firsthomo}. This latter kernel is easily identified
with $\oplus_e H_1(\widehat{Z}(e)) \otimes \scrA_{\ell^*(e)}$ which
is naturally identified with
$\bigoplus_{e;\ell\in\stackrel{\circ}{e}} H_3(R(\ell))$.
\end{proof}

Now, let's compute the image of $\rho_* \colon H_3(\tY) \to
H_3(Y)$.  We have the long exact commutative diagram associated to
the sequences
\begin{eqnarray*}
\tY_2 \subset \tY_3 \to (\tY_3, \tY_2) & {\simeq} &
(\widetilde{\nu}_3, \pa \widetilde{\nu}_3) \\
Y_2 \subset Y_3 \to (Y_3, Y_2) & {\simeq} & (\nu_3, \pa {\nu}_3).
\end{eqnarray*}
We obtain
$$\begin{array}{ccccccccc}
0 & & & & & & & \\
\downarrow & & & & & & & \\
 \oplus_e H_1(Z_e) \otimes \scrA_{\ell^*(e)} & &  & & & & & \\
 \downarrow                                  & &  & & & & & \\
H_{3}(\tY_2) & \to & H_3(\tY) & \to & 0 & & & \\
\downarrow & & \downarrow & &  \downarrow & & & \\
H_3(Y_2) & \to & H_3(Y) & \to
& H_3(\nu_3, \pa {\nu}_3) & \to & H_2(Y_2) & \to & H_2(Y) \\
 \downarrow    &  &        &     & & & & \\
0    & &   &     & & & &
\end{array}$$
yielding
$$\Image(H_3(\tY) \to H_3(Y)) = \Image(H_3(Y_2) \to H_3(Y)) .$$
\end{proof}

It is not true in general that $\rho_* \colon H_3(\tY)\to H_3(Y)$ is onto, or equivalently
that $H_3(Y_2)\to H_3(Y)$ is onto. From the exact sequence of the pair
$(Y,Y_2)$ and the fact that $(Y,Y_2)$ is excisively equivalent to
$(\nu_3,\partial \nu_3)$ we see that the cokernel of $H_3(Y_2)\to H_3(Y)$ is
identified with the kernel of $H_3(\nu_3,\partial \nu_3)\to H_2(Y_2)$. Of
course, since each component of $\nu_3$ is a cone, this kernel is identified
with the kernel of the map induced by the inclusion $H_2(\partial \nu_3)\to
H_2(Y_2)$. Let's compute rationally. Arguing as in the proof of
Lemma~\ref{2.22} using $Y_0\subset Y_1\subset Y_2$ we see that $H_1(Y_2;\Q)=0$
and
$$H_2(Y_2;\Q)=\Ker \left(\oplus_{v\in V}(\Delta)\Zee\langle v\rangle\to
N\otimes \Q\right).$$ Clearly, this map is onto rationally and hence the rank
of $H_2(Y_2)$ is $\#V(\Delta)-n$. For each two-face $f$ of $\partial \Delta$,
the boundary of the regular neighborhood $\partial \nu(f)$ is homotopy
equivalent to $V(f)^*$, the complement of the fixed point $0_f$ in $V(f)$.
Similar arguments to the ones for $Y_2$ show that the rank of $H_2(V(f)^*,\Q)$
is given by the number of vertices of $f$ minus $3$. (Indeed the map between
these groups is easily determined from the combinatorial configuration of
vertices and two-faces.) In any event if the sum over the two-faces of the
number of vertices of the face minus three is greater than the number of
vertices of $\Delta$ minus $n$, then the kernel of this map is non-trivial and
hence $H_3(\tY)\to H_3(Y)$ will not be onto. As an explicit example where this
map is not onto, we have the $4$-cube which has $24$ two-faces, each with four
vertices. For each two-face $f$ the affine three-fold $V(f)$ is a (complex)
cone on the quotient of $\Pee^1\times \Pee^1$ by an involution preserving the
factors. Thus, $V(f)^*$  has second homology of rank $1$. The cube itself has
$16$ vertices so that $H_2(Y_2)$ has rank $12$. Hence, in this case the kernel
of the map has rank at least $12$, and in fact has rank $13$.

\begin{corollary}\label{kernelformula}
In the hypersurface case the rank of the kernel of $\rho_* \colon H_3(\tY)\to
H_3(Y)$ is
$$2\sum_{e\in E}\ell^*(e)\ell^*(e^\vee).$$
For complete intersections of ample divisors the rank of the kernel
of the map $\rho_* \colon H_3(\tY)\to H_3(Y)$ is given by
$$2\sum_{e\in E}\ell^*(e)\left( \sum_{J \subset I;\ J\not=\emptyset}
(-1)^{3 - |J|} \ell^* (\sum_{j\in J}e^\vee_j)\right).$$ In both
cases the Hodge structure on this kernel is of type $(2,1)$ and
$(1,2)$ each piece being of half the rank.
\end{corollary}

\begin{proof}
According to Corollary~\ref{dim1&2hyper} in the hypersurface case
the rank of $H_1(\widehat Z(e))$ is $2\ell^*(e^\vee)$. According to
Corollary~\ref{dim1&2ci} the formula for the rank of $H_1(\widehat
Z(e))$ in the general case is
$$2 \sum_{J \subset I;\ J\not=\emptyset} (-1)^{3 - |J|}
\ell^*(\sum_{j \in J} e^\vee_j).$$  The result is immediate from
this.
\end{proof}

\subsection{The tangent space to the space of polynomial deformations}

\begin{lemma}
Suppose $M$ is a smooth $n$-dimensional Calabi-Yau manifold and
let $D \subset M$ be a smooth divisor. Suppose the pair $(M,D)$
deforms.  Let
$$\alpha_{M} \in H^1(M, T_M) = H^{n-1,1}(M)$$
be the Kodaira-Spencer class of the deformation of $M$.  Then the
restriction of $\alpha_{M}$ to a class in $H^{n-1,1}(D)$ vanishes.
\end{lemma}

\begin{proof}
Suppose we have a deformation of the pair $(M,D)$.  We have the
following sequence in cohomology
$$\begin{array}{rcl} H^0(D, \nu_{D \subset M}) \to H^1(D, T_D)
\to & H^1(D, T_{M}|_D) & \to H^1(D, \nu_{D \subset {M}})
\\ & \uparrow & \\ & H^1(M, T_{M}) &
\end{array}$$
where, via the Kodaira-Spencer mapping, $H^1(M, T_{M})$
corresponds to the tangent space to deformations of $M$ and
$H^1(D, T_D)$ to those of $D$. Deformations of the pair $(M,D)$
are given by a pair of classes $\alpha_{M} \in H^1(M,T_{M})$ and
$\alpha_D \in H^1(D,T_D)$ with the same image $\alpha' \in H^1(D,
T_{M}|_D)$. Thus, the composite mapping
$$H^1(M, T_{M}) \to H^1(D, T_{M}|_D)$$
sends $\alpha_M$ to zero.  Moreover, we have the sequence
$$\begin{array}{ccc}
H^1(M, \Omega_{M}^{n-1}) & \to H^1(D, \Omega_{M}^{n-1}|_D) \to &
H^1(D, \Omega_D^{n-1})
 \\
|| & & || \\
H^1(M, T_{M}) & & H^1(D, K_D)
\end{array}$$
so that the composition map $$H^1(M,T_{M}) \to H^1(D,T_{M}|_D) \to
H^1(D, T_D)$$ is identified with the natural restriction map
$$H^{n-1,1}(M) \to H^{n-1,1}(D)$$
This shows that if the pair $(M, D)$ deforms, then the restriction
map on $H^{n-1,1}$ is trivial. Conversely, if the image of
$\alpha_{M}$ in $H^1(D, T_{M}|_D)$ goes to zero in $H^1(D, \nu_{D
\subset M})$ then there is a class $\alpha_D \in H^1(D,T_D)$ with
the same image as $\alpha_{M}$ in $H^1(D,T_{M}|_D)$.  The pair
$(\alpha_{M},\alpha_D)$ then gives an infinitesimal deformation of
$(M,D)$.
\end{proof}

Now let us apply this to the situation of Calabi-Yau threefolds in
toric varieties.

\begin{corollary}
Let $\alpha\in H^1(\tY,T\tY)=H^{2,1}(\tY)$ be the Kodaira-Spencer
class of a deformation of $\tY$ induced by taking the preimage
under the map $\rho\colon \Pee_\scrT\to \Pee_\Delta$ of a
polynomial deformation of $Y\subset X$ (i.e., a deformation
obtained by varying the coefficients of the polynomials cutting
out $Y$ as a complete intersection). Then for each edge $e$ of
$\pa \Delta$ and any lattice point $\ell\in \stackrel{\circ}{e}$,
the element $\alpha\in H^{2,1}(\tY)$ restricts to zero in
$H^{2,1}( R(\ell))$.
\end{corollary}

\begin{proof}
Since the deformation $\{Y_t\}$ of $Y$ takes place in
$\Pee_\Delta$ we have the family of curves $\widehat{Z}_t(e) =
Y_t\cap\overline \scrO_e$ and their preimages which are surfaces
$R_t(\ell)\subset \tY_t$ deforming $R(\ell)$. Thus, the given
deformation of $\tY$ lifts to a deformation of the pair $(\tY,
R(\ell))$. Applying the previous  lemma gives the result.
\end{proof}

\begin{corollary} Let
$H^{2,1}_{\rm poly}(\tY)$ be the tangent space of the polynomial
deformations of $\tY$ inside $H^{2,1}(\tY)$. This subspace is
identified with the kernel of the map
$$H^{2,1}(\tY)\to \bigoplus_{\ell\in N\cap \pa\Delta^{(1)}\setminus V(\Delta)}
H^{2,1}(R(\ell)).$$
\end{corollary}

\begin{proof}
The previous result shows that $H^{2,1}_{\rm poly}(\tY)$ is
contained in the kernel of the restriction mapping. In
Corollary~\ref{kernelformula} we showed that the rank of the kernel
of the restriction mapping is equal to the rank of the space of
non-polynomial deformations (i.e., the correction term) as computed
%[[specific reference]]
by \cite{BB}.
\end{proof}

%[[$W_{-3}$ Hodge structure comments]]
%
%[[restriction to ruled surfaces to give Coker]]
%
%****************
%
%Batyrev-Borisov hypersurface formulas:
%
%$$h^{1,1} = \ell(\Delta) - 5 - \left( \sum_{ \{\Theta
%\subset \Delta \, | \, \codim \Theta = 1\} } \ell^*(\Theta)
%\right) + \left( \sum_{ \{\Theta \subset \Delta \, | \, \codim
%\Theta = 2\} } \ell^*(\Theta) \cdot \ell^*(\Theta^\vee) \right)$$
%$$h^{2,1} = \ell(\Delta^\circ) - 5 - \left( \sum_{ \{\Theta \subset
%\Delta \, | \, \codim \Theta = 4\} } \ell^*(\Theta^\vee) \right) +
%\left( \sum_{ \{\Theta \subset \Delta \, | \, \codim \Theta = 3\}
%} \ell^*(\Theta) \cdot \ell^*(\Theta^\vee) \right)$$
%
%*********************************************
%
%Batyrev-Borisov complete intersection formulas ($d-r=3$): Complete
%intersection of pullbacks of ample divisors; $\Delta$ is
%combinatorially dual to each of $\Delta^\circ, \Delta_1^\circ,
%\ldots, \Delta_r^\circ$.
%
%$$h^{1,1} = \# V + \sum_{e \in E} \ell^*(e) + \sum_{f \in F} \ell^*(f) -
%d + \sum_{f \in F} \ell^*(f) \left( \sum_{J \subset I} (-1)^{r -
%|J|} \ell^*( \sum_{j \in J} f_j^\vee ) \right)$$
%
%$$h^{2,1} = \sum_{i=1}^r \left( \sum_{J \subset I} (-1)^{r - |J|}
%\ell^*(\Delta_i^\circ + \sum_{j \in J} \Delta_j^\circ) \right) - d -
%\sum_{v \in V} \left( \sum_{J \subset I} (-1)^{r - |J|}
%\ell^*(\sum_{j \in J} v_j^\vee) \right) $$ $$ \ \ + \sum_{e \in E}
%\ell^*(e) \left( \sum_{J \subset I} (-1)^{r - |J|} \ell^*(\sum_{j
%\in J} e^\vee_j) \right) $$
%
%
%
%\section{Cohomology and Hodge structure of complete intersection
%Calabi-Yau threefolds}

\section{$K^0$ and $K^1$ of Calabi-Yau threefolds in terms of homology}

The 7-dimensional stage of the Postnikov tower for $BSU$, denoted
$BSU^{(7)}$, contains the homotopy groups $\pi_i(BSU)$ for $i \leq
7$ and has trivial groups in degrees $\geq 8$. In particular,
$$\pi_i(BSU^{(7)}) \simeq \left\{ \begin{array}{cl} \Zee & i = 4,
6 \\ 0 & \mbox{otherwise} . \end{array} \right.$$  If $X$ is a CW complex of
dimension $\leq 7$, then the natural map from
$$[X, BSU] \to [X, BSU^{(7)}]$$
is a bijection.

$BSU^{(7)}$ is given by a fibration \begin{equation}
\label{eqnKBKfibrn} K(\Zee,6) \to BSU^{(7)} \to K(\Zee,4)
\end{equation} with $k$-invariant $\delta Sq^2 \iota_4$ the unique nontrivial
element (of order two) in $H^7(K(\Zee,4);\Zee)$. Thus we also have
a fibration
$$BSU^{(7)} \to K(\Zee,4) \times K(\Zee,6) \to K(\Zee/2,6)$$
where the first map is given by the second and third Chern classes
$(c_2, c_3)$ and the second map is given by $Sq^2 \iota_4 +
[\iota_6]_2$, with $\iota_4$ and $\iota_6$ being the fundamental
classes of the factors and $[ \ \cdot \ ]_2$ denoting reduction
mod 2.

This means that for any CW-complex $X$ of dimension $\leq 7$ we
have that
$$[X, BSU] = [X, BSU^{(7)}]$$
sits in an exact sequence
$$\hspace{-.2in} H^3(X;\Zee) \oplus H^5(X;\Zee) \stackrel{\alpha}{\to}
H^5(X;\Zee/2) \to [X, BSU] \to H^4(X;\Zee) \oplus H^6(X;\Zee)
\stackrel{\beta}{\to} H^6(X; \Zee/2)$$ where $\alpha(a,b) = Sq^2 a
+ [b]_2$ and $\beta(c,d) = Sq^2 c + [d]_2$.  In the special case
that $M$ is a closed, oriented 6-manifold, $H^6(M;\Zee)$ is
torsion-free, and hence the mod 2 reduction $H^5(M;\Zee) \to
H^5(M;\Zee/2)$ is onto.  Thus, we have
$$0 \to [M, BSU] \stackrel{(c_2,c_3)}{\to} H^4(M;\Zee) \oplus
H^6(M;\Zee) \to H^6(M; \Zee/2) \to 0 .$$ This proves

\begin{lemma}
Let $M$ be a closed oriented 6-manifold. Then $[M,BSU]$ is isomorphic to
$$\{(c_2,c_3) \in H^4(M;\Zee) \oplus
H^6(M;\Zee) \vert Sq^2 c_2 = [c_3]_2 \}$$ where the isomorphism is given by
taking the  2nd and 3rd Chern class.
\end{lemma}

An examination of low dimensional examples allows one to extend this result to
$\widetilde {K^0}(M)= [M,BU]$:

\begin{lemma}\label{c1c2c3}
Let $M$ be a closed oriented 6-manifold.  Then $\widetilde {K^0}(M)$ is
isomorphic to
$$\{(c_1,c_2,c_3) \in H^2(M;\Zee) \oplus H^4(M;\Zee) \oplus
H^6(M;\Zee) \vert Sq^2 c_2 = [c_3]_2 +c_1c_2+c_1^3\}$$ where the isomorphism is
given by taking the 1st, 2nd, and 3rd Chern class.
\end{lemma}

\begin{corollary}
If $M$ is a Calabi-Yau threefold, then $$\Tor K^0(M) = \Tor
H^2(M;\Zee) \oplus \Tor H^4(M;\Zee)$$ and
$$\rk  K^0(M) = \sum_{* \leq 3} \rk \left( \left. H^{2*}(M; \Zee)
\right/ \mbox{Torsion} \right) \ .$$
\end{corollary}
\begin{proof}
The rank statement is well-known.  The torsion statement follows
from the above since
$$Sq^2 \colon H^4(M;\Zee) \to H^6(M;\Zee)$$
is zero.  This follows from the fact that this map is the cup product with
$w_2(M)$ which, for any almost complex manifold, is the reduction mod 2 of
$c_1(M)$. The Calabi-Yau condition means that $c_1(M)=0$.
\end{proof}
In other words, we have that, for a Calabi-Yau threefold $M$, the
even K-group is
$$K^0(M) \simeq \Zee \oplus H^2(M;\Zee) \oplus H^4(M;\Zee) \oplus
2 \cdot H^6(M;\Zee) .$$

Similar arguments show that \begin{lemma} Let $M$ be a closed
oriented 6-manifold.  Then $K^1(M) = \widetilde{K^0}(\Sigma M)$ is
$$H^1(M;\Zee) \times [\Sigma M, BSU^{(7)}]$$ and we have an exact
sequence
$$0 \to H^5(M;\Zee) \to [\Sigma M, BSU^{(7)}] \to H^3(M;\Zee) \to 0 $$
coming from the fibration (\ref{eqnKBKfibrn}). \end{lemma}
\begin{corollary}
When $M$ is a Calabi-Yau threefold the extension class is trivial.
\end{corollary}
\begin{proof}
Once again
$$Sq^2 \colon H^3(M;\Zee) \to H^5(M;\Zee/2)$$
is the reduction mod 2 of cupping with $c_1$, which is consequently trivial.
Thus,  by Lemma~\ref{c1c2c3}, for any class $a\in H^3(M;\Zee)$ there is a
unique (up to isomorphism) bundle over $\Sigma M$ with $c_2=a$ and $c_1=c_3=0$.
The Whitney sum formula shows that these bundles form a subgroup in $\widetilde
{K^0}(\Sigma M)$, and hence this construction gives a splitting of the
sequence.
\end{proof}
In other words, we have that, for a Calabi-Yau threefold $M$, the
odd K-group is
$$K^1(M) \simeq H^1(M;\Zee) \oplus H^3(M;\Zee) \oplus H^5(M;\Zee) .$$

\end{document}